\documentclass{article}
\pdfoutput=1
\usepackage{lipsum}
\usepackage{amsfonts}
\usepackage{graphicx}
\usepackage{epstopdf}
\usepackage{algorithmic}
\usepackage{mathtools}
\usepackage{todonotes}
\usepackage{subcaption}
\usepackage{hyperref}
\ifpdf
  \DeclareGraphicsExtensions{.eps,.pdf,.png,.jpg}
\else
  \DeclareGraphicsExtensio\pdfoutput=1ns{.eps}
\fi
\newcommand{\beqn}{\begin{eqnarray}} 
\newcommand{\eeqn}{\end{eqnarray}}

\newcommand{\MSR}{MS\&R\ }

\newcommand{\TheTitle}{The shifted proper orthogonal decomposition:\\
A mode decomposition for multiple transport phenomena} 
\newcommand{\TheAuthors}{J. Reiss, P. Schulze, and J. Sesterhenn}

\title{{\TheTitle}}

\author{
  J{.} Reiss\thanks{Institut f\"ur Str\"omungsmechanik und Technische Akustik, TU Berlin, Germany, \texttt{$\{$reiss,sesterhenn$\}$@tnt.tu-berlin.de}.}
  \and
  P{.} Schulze\thanks{Institut f\"ur Mathematik, TU Berlin, Germany, \texttt{$\{$pschulze,mehrmann$\}$@math.tu-berlin.de}.}
  \and
  J{.} Sesterhenn\footnotemark[1]
  \and
  V{.} Mehrmann\footnotemark[2]
}

\usepackage{amsopn}

\ifpdf
\hypersetup{
  pdftitle={\TheTitle},
  pdfauthor={\TheAuthors}
}
\fi

\begin{document}

\maketitle

\begin{abstract}
\noindent Transport-dominated phenomena 
provide a challenge for common mode-based model reduction approaches. 
We present a model reduction method, which is suited for these kinds of systems. It extends the proper orthogonal decomposition (POD) by introducing  time-dependent shifts of the snapshot matrix. 
The approach, called shifted proper orthogonal decomposition (sPOD), features a determination of the {\it multiple} transport velocities and a separation of these. 
One- and two-dimensional test examples reveal the good performance of the sPOD for transport-dominated phenomena and its superiority in comparison to the POD.
\end{abstract}


  
\section{Introduction\label{sec:intro}}
Model reduction is an important method to derive low-di\-men\-si\-onal models from experimental or numerical data.
These reduced-order models often allow fast simulations commonly used for control, optimization, and parameter studies and are nowadays a necessary tool in many fields.
Further, these models yield a better understanding of the dynamical process by identifying the essential dynamics.
Formally, the goal is to obtain a low-dimensional description, which approximates the mapping from a set of inputs to a set of outputs.
Among others, inputs can be design parameters, system conditions, or controls.
Common outputs are performance or physical quantities like measurements or even full flow solutions.

A major class of reduced-order modeling approaches is given by input-output interpolation methods, which do not aim in reducing the internal dynamics but which are only based on the input-output behavior, usually described by a transfer function (see \cite{AntBG10,BeaG09,GugAB06,MayA07} and references therein).
 These methods are successfully applied to linear systems, but extensions to nonlinear systems are rare and have some drawbacks as low computational efficiency (cf. \cite{FujS10}).

In contrast, the most common model reduction techniques for nonlinear systems are based on  a superposition of modes describing the system state.
Examples are reduced basis methods \cite{DihH15,Eft11,Gre12,RozHP08,UrbP12}, balanced proper orthogonal decomposition  (POD) \cite{TuR12,WilP02}, dynamic mode decomposition \cite{CheTR12,SchmidSesterhenn2008,Tu13}, and POD \cite{HinV05,RatP03}, which technically reduces to a singular value decomposition (SVD).
The SVD provides the best low-rank approximation of a matrix with respect to the $2$-norm. 
Depending on the application, the POD has several synonyms, as for instance, principal component analysis \cite{Jol86} or Karhunen-Loeve decomposition \cite{DehM08}.
The number of modes needed to obtain a sufficiently accurate approximation of the system in question is crucial for the computational cost of evaluating the reduced-order model (online stage) and accordingly for its usefulness.
Furthermore, nonlinear reduced-order models are often combined with hyperreduction methods, e.\,g., EIM \cite{BarMNP04} or DEIM \cite{ChaS10}, to achieve an efficient offline/online decomposition \cite{HaaO11}.

Transport-dominated phenomena are usually a challenge for mode-based methods, since their dynamical behavior cannot be captured accurately  by a few spatial modes in a dyadic structure (cf. \cite{ConI12,Luc01}).
Recently, there have been some efforts in inventing model reduction methods suitable for the efficient description of transport-dominated phenomena.
Usually some time-dependent shift is introduced to compensate the transport.
It is used in the framework of symmetry reduction (cf. \cite{BeyT04,FedAR15,OhlR13,RowKML03}), where the translation is accounted for by applying a Lie group action to a symmetry-reduced or frozen solution.
In \cite{OhlR13}, for instance, the framework has been analyzed for nonlinear parameter-dependent evolution equations and applied to a numerical simulation of the Burgers equation.
A different approach is presented in \cite{IolL14} where the transport is incorporated by using a coordinate mapping which is related to the solution of Monge-Kantorovich optimal mass transport problems.
The approach is illustrated by means of snapshots of shallow water waves and of a hurricane.
A methodology based on $L^1$ norm minimization has been applied in \cite{AbgAC16} which shows much better results for hyperbolic problems than the commonly applied $L_2$ norm minimization of the error.
The minimization is based on a set of dictionaries which are computed in an offline phase.
There are various other approaches which aim to efficiently reduce transport phenomena, see for instance \cite{CagMS16,Car15,GerL14,MojB17,SesterhennShahirpour2016,TadPQ15,Wel15}.

Multiple transport velocities are less studied. Just one shift as it is used in most of the cited works is not sufficient if different transport velocities are present, as it is common in technical applications.
An efficient and general model reduction methodology for {\it multiple} transport phenomena is still missing.
With this contribution we aim at improving this situation  by introducing the shifted proper orthogonal decomposition (sPOD).

The key tool of the new approach is a time-dependent shift in combination with a procedure to {\it separate} different transports within the system.
The dominant transport velocities are determined by front tracking or by considering the dependence of the singular values of the shifted snapshot matrix on the time-dependent shift. This shift structurally extends the dyadic structure of POD to systems with transport and allows thereby a better approximation.
The method is purely data based in contrast to the symmetry reduction methods, offering a wide range of applicability.

In this paper, we focus on obtaining a very low-dimensional representation of the solution.
The corresponding low-dimensional subspace can be used to build a reduced-order model (ROM), for instance by an interpolation procedure as in \cite{LemMRMS15} or by a Galerkin projection of the original model. However, the construction of a reduced-order model is not within the scope of this paper.

During the review process we became aware of the work \cite{RimMoeLeVeque2017}, which is similar to this work using SVDs in different velocity frames.
The authors of \cite{RimMoeLeVeque2017} utilize a greedy approach in order to decompose snapshot matrices with multiple transport velocities. They consider the linear wave equation with two transported quantities and their method was shown to yield a decomposition with just a few modes outperforming the POD. However, the greedy approach is not able to describe this linear transport with the minimum number of modes indicated by the analytic solution. This is because structures which are extracted by the first greedy iteration cannot be re-attributed to a different velocity frame, which is necessary to obtain low rank approximations for hyperbolic cases with multiple transports. In contrast, this flexibility is given by the new sPOD algorithm which is able to decompose the same example with the minimum number of modes and thus finds the analytic solution up to any given achievable tolerance (cf. section \ref{sec:multipleTransports}).

In the following section, we motivate and derive the method by means of examples in a one-dimensional domain, which are a challenge for common model reduction techniques as the POD. For this, first we consider one moving signal and show how to reduce it. Then, a more complex system with two different transports is considered and a procedure is derived to separate the different velocity components.
Then the  detection of the velocities is discussed and a case of crossing shocks is investigated.
In section \ref{sec_twoD} we apply the new method to a two-dimensional test case from fluid mechanics with transported developing vortex pairs with non-trivial velocities.
Finally, we conclude and give an outlook on ongoing and future work.

\section{One-dimensional model problems}
\label{sec_oneD}

In this section, first the idea of the sPOD approach will be developed considering a one-dimensional example.
As a test problem, the linear wave equation
\begin{equation}
\begin{aligned}
\partial_t \rho  +  \rho_0   \partial_x u  &=& 0, \\
\partial_t  u   +  c^2/\rho_0 \partial_x \rho  &=& 0 ,
\end{aligned}
\label{linearWave}
\end{equation}
is considered on a periodic domain $x= (0,L] $.
Here, $u$ is the velocity, $\rho$ is the density (fluctuation), $\rho_0$ is a reference density, and $c$ the speed of sound.
The general solution can be written as
\beqn
q(x,t) = \begin{bmatrix}
\rho\left(x,t\right)\\
u\left(x,t\right)
\end{bmatrix}
=
q_+ (x - c t)
\begin{bmatrix}
\rho_0\\
c
\end{bmatrix}
+
q_- (x  - (-c) t)
\begin{bmatrix}
 \rho_0\\
- c
\end{bmatrix}
\label{solRiemann}
\eeqn
with arbitrary initial conditions $q_\pm (x) $  for the  two transported quantities, the Riemann invariants.
The corresponding transport velocities are $ c_{\pm} = \pm c$.
By choosing
\beqn	
 q_\pm(x-(\pm ct)) = \frac 1 2  \sum_{n=0}^\infty \beta_n \cos( k_n(x-(\pm ct)) +\Theta_n  ), \label{waveRI}
 \eeqn
 with $k_n= n\, 2\pi/L$, the solution can be rewritten  as
\begin{eqnarray}
q(x,t) = 
\sum_{n=0}^\infty \beta_n 
\begin{bmatrix}
  \rho_0\cos(k_n\,   ct )   \cos(k_n\,  x +\Theta_n) 
  \\
 c\sin(k_n\,  ct )   \sin(k_n\,  x +\Theta_n) 
\end{bmatrix}
\label{solString}
%
\end{eqnarray}
with constants $\beta_n$ and $\Theta_n$. Every non-zero amplitude $\beta_n$  yields  a mode of a vibrating string.
For the remainder of this paper, we set $\rho_0 = 1$, $c= 1$, and $L=1$.
\medskip

\noindent As stated in section \ref{sec:intro}, model reduction methods often build on describing a dynamical system by a superposition of a small number of modes. One of the most popular approaches is the POD, which aims at approximating the solution by a linear combination of orthonormal
modes $\phi_l$
\beqn
q(x_i, t_n ) \approx \sum_l \alpha_l(t_n)\phi_l(x_i)  \label{pod}
\eeqn
with time-dependent coefficients $\alpha_l$. Usually, a snapshot matrix $X_{ij} = q(x_i,t_j)$ is introduced, which is a space-time-discretized solution.
A low-dimensional representation minimizing the approximation error in the $2$-norm is determined via a singular value decomposition (SVD) of the snapshot matrix $X$. Often the numerical effort of the SVD is decreased by reducing the set of discrete time and spatial points  \cite{BarAA05} in a way that does not significantly effect the quality of the modes, i.\,e., the dynamical behavior of the system still needs to be captured. 
The squared singular values determine the mean square amplitude of the corresponding modes, which coincides with the kinetic energy if the mode represents the velocity \cite{LucNS09}. Consequently, strongly decaying singular values allow a good representation with a few modes.

The ability of the POD to describe the solution of (\ref{linearWave}) eminently depends on the initial condition.
If a vibrating string with a dominant frequency and a few harmonics is chosen, i.e., as in \eqref{solString} with a few dominant $\beta_n$,  the POD will find these few modes and, consequently, it will deliver an accurate low-dimensional representation.
However, if a transported quantity with high gradient is given, i.e., as in \eqref{solRiemann} with strongly localized initial conditions $q_\pm$, the singular value decay is rather gradual and one is forced to use a high number of modes to get a reasonable description.
 This variability makes the example well-suited to develop a method, which can handle transported quantities while including the classical POD as a special case.

\subsection{One transported quantity\label{sec:oneTransport}}

\begin{figure}
\begin{center}
\includegraphics [width=0.49\linewidth,clip=true ,  viewport = 0 5 535 420]{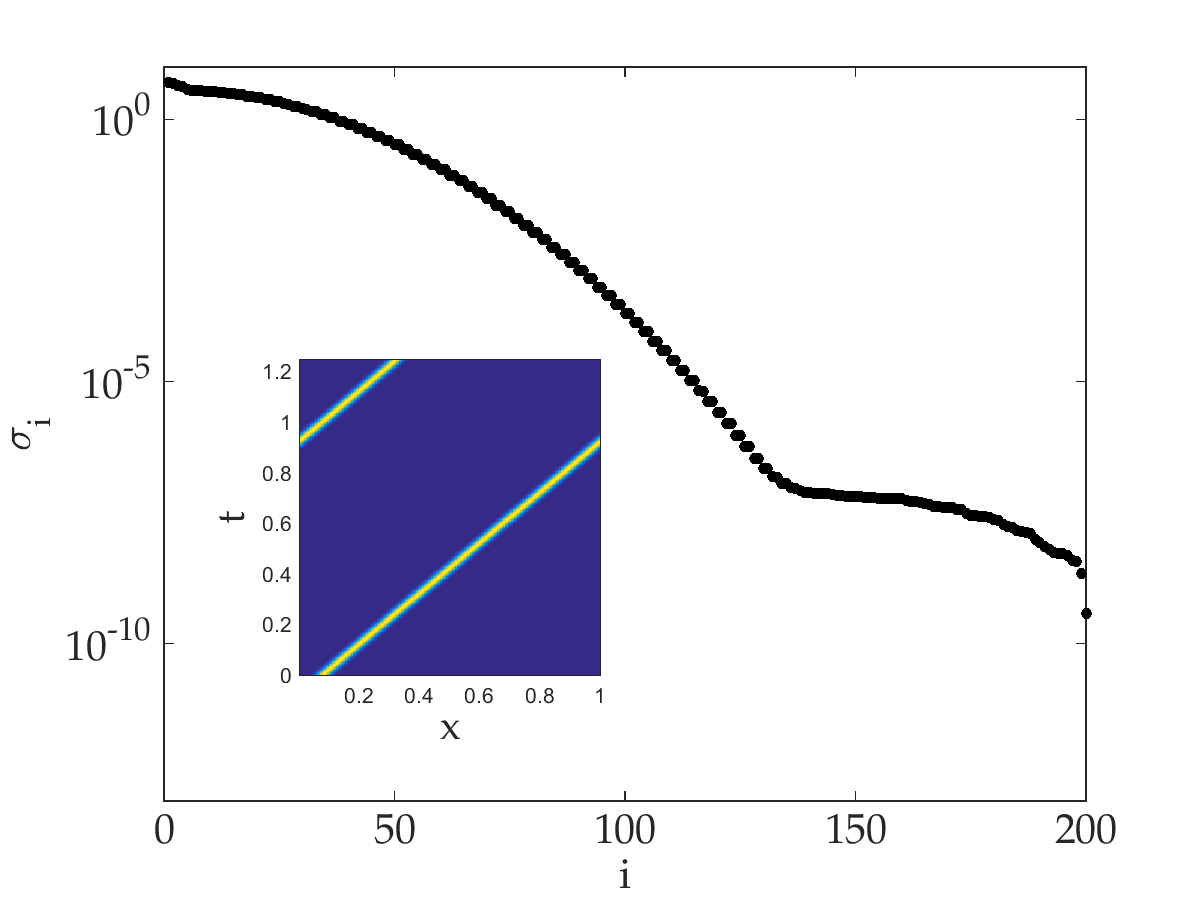}
\includegraphics [width=0.49\linewidth,clip=true ,  viewport = 0 5 535 420]{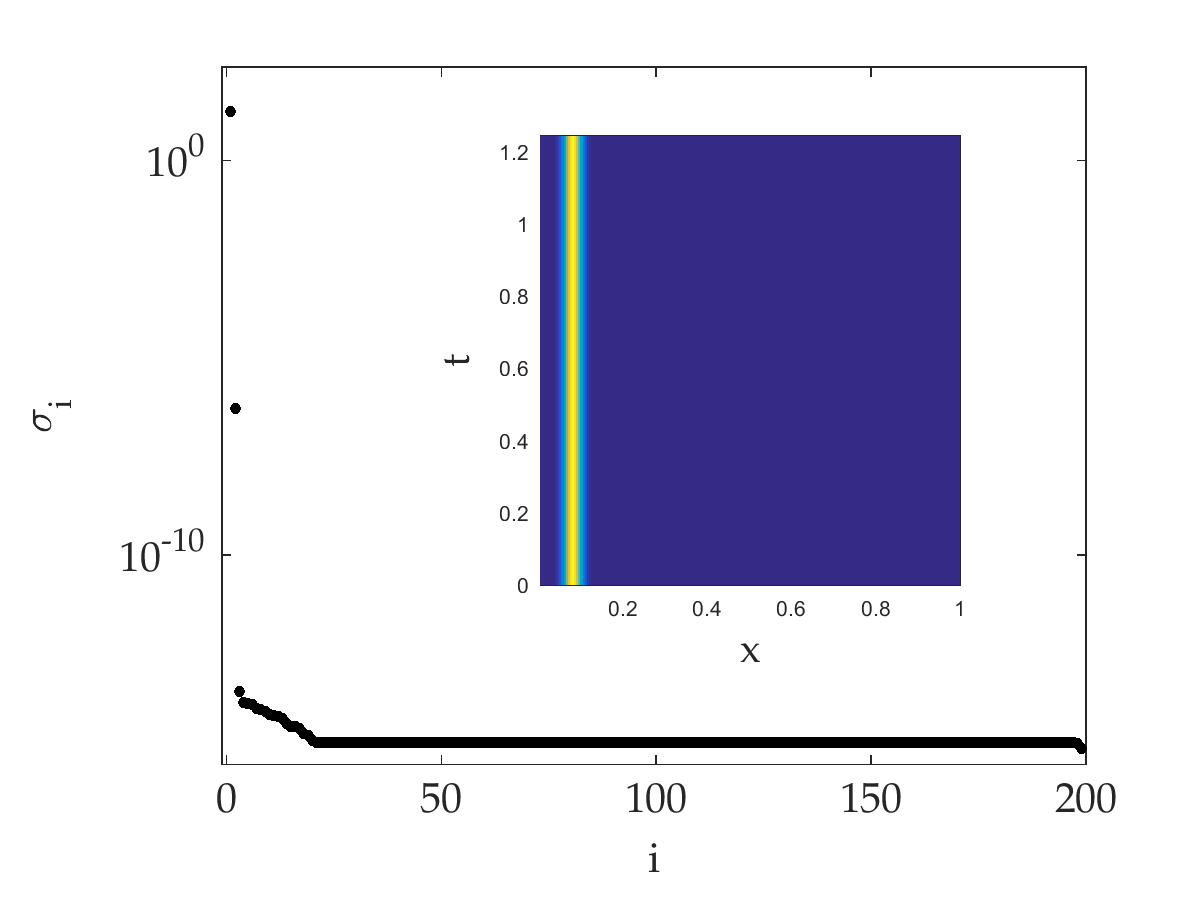}
\end{center}
\caption{
Single transported quantity. Left: Singular value decay (outer picture) and density in the space time diagram (inner picture). 
The spectrum shows a very slow decay.
Right: If the solution is described  in a  co-moving frame, i.e., the solution is shifted at every time to compensate the transport, the singular spectrum shows a rapid decay. 
One mode is sufficient to describe the flow (apart from a small numerical error of the shift procedure).
}
\label{fig:rightGoing_Spectrum}
\end{figure}

Setting $q_-(x) =0$  in  the  solution \eqref{solRiemann} of the linear wave equation \eqref{linearWave} leads to a single transported quantity,
which is given by  a shift of the initial condition $q_+$ with the transport velocity $c$. To provide a challenging case for POD, a sharp Gaussian pulse of $q_+(x ) = \exp ( - (x-x_0)^2/\delta^2 )  $ is chosen with $\delta = L/50$. The analytic solution is shown in Fig.\ \ref{fig:rightGoing_Spectrum}, left, for the time interval $[0, 1.25 L/c] $. For the discrete snapshot matrix $X_{ij}  = q(x_i,t_j ) $, we chose 250 equidistant time and 200 space points.
In this subsection we only consider the density, if not stated otherwise.

For the considered example, a very slow decay of the singular values is observed (cf.\ Fig.\ \ref{fig:rightGoing_Spectrum}, left). 
Even though the analytic solution can be formulated by only one transported Gaussian pulse, many POD modes are needed for a good representation. The reason for this is that the structure
(\ref{pod}) is not adequate to describe a solution of the form (\ref{solRiemann}).
Roughly speaking the dyadic products  describe rectangular (with respect to the space time diagram) structures well, while diagonal structures such as in Fig.\ \ref{fig:rightGoing_Spectrum}, left, cannot be well represented as a dyadic product. This explains the common failure of the mode based model reduction for transport phenomena (cf.\ \cite{CagMS16,RimMoeLeVeque2017}).
However, if the velocity is known or can easily be detected, the solution is near at hand.
A time-dependent shift which compensates the transport velocity yields a structure that can be written as a single dyadic product, which in turn corresponds to the description by just one mode,  Fig.  \ref{fig:rightGoing_Spectrum}, right.

To this end, we introduce the continuous shift operator $\mathcal{T}^{c} f\left(x,t\right) \vcentcolon = f\left(x - ct,t\right)$ acting on a function $f\left(x,t\right)$. The discrete analogue of $\mathcal{T}^{c}$ is denoted with $T^c\left(\cdot\right)$ which acts on a snapshot matrix of a time- and space-dependent function and shifts the $i$th column by $ct_i$ in space.
Thus, for one transport we  seek approximations
\beqn
q(x_i,t_n) \approx T^{c}\left(\sum_l \alpha_l(t_n)\phi_l(x_i) \right). \label{oneAnsatz}
\eeqn
Clearly such an approximation can be constructed by applying the inverse shift $T^{-c} $ into the co-moving frame and performing an SVD.
 For our example of the single transported quantity, we obtain,  $T^{-c}\left( q(x,t)\right) = q( x  ,0 )$ which is described by one mode.
 Note that this shift demands an interpolation in the general case where $ct_i$ is not an integer multiple of the grid size.
Note further that the approach also works with a non-constant transport velocity, since the shift can be set for each time step separately.

In a similar manner, shifts have been used in other reduction frameworks, e.\,g., for the model reduction of a combustion \cite{LemMRMS15} and in a symmetry reduction framework \cite{BeyT04,FedAR15,RowKML03}.
The latter one treats partial differential equations (PDEs) which are equivariant with respect to a group action which in turn induces symmetries in the solution space.
Symmetry-reduced surrogate models on the PDE level are obtained by exploiting the equivariance of the original PDE.
In contrast, we are looking for a general decomposition of a snapshot matrix by shifted modes without requiring equivariance of the original PDE or symmetries in their solutions.

\subsection{Multiple transported quantities\label{sec:multipleTransports}}

Many relevant  systems 
feature multiple transported quantities.
We extend \eqref{oneAnsatz} to a \textit{multi-frame decomposition} of the form
\beqn
q(x_i,t_n) \approx \sum_{k=1}^{N_{s}}T^{c_k}\left(\sum_l \alpha_l^k(t_n)\phi_l^k(x_i) \right) \label{multiAnsatz}
\label{eqn:multi_struc}
\eeqn
where $N_{s}$ is the number of transport velocities.
We seek to decompose a given  space-time field into this structure.
To this end, a best fit of (\ref{eqn:multi_struc}) in the vector 2-norm  will be used, where the ansatz modes will be constructed  by shifting and reducing the field data to identify low rank structures  in each co-moving frame.
The reduction and best fit is iterated to decrease the cross-influence between the different transports and to obtain a clear separation.

\begin{figure}[!htb]
\centering
\includegraphics[scale = 0.6,clip=true,trim = 0.5cm 0cm 1.5cm 0cm]{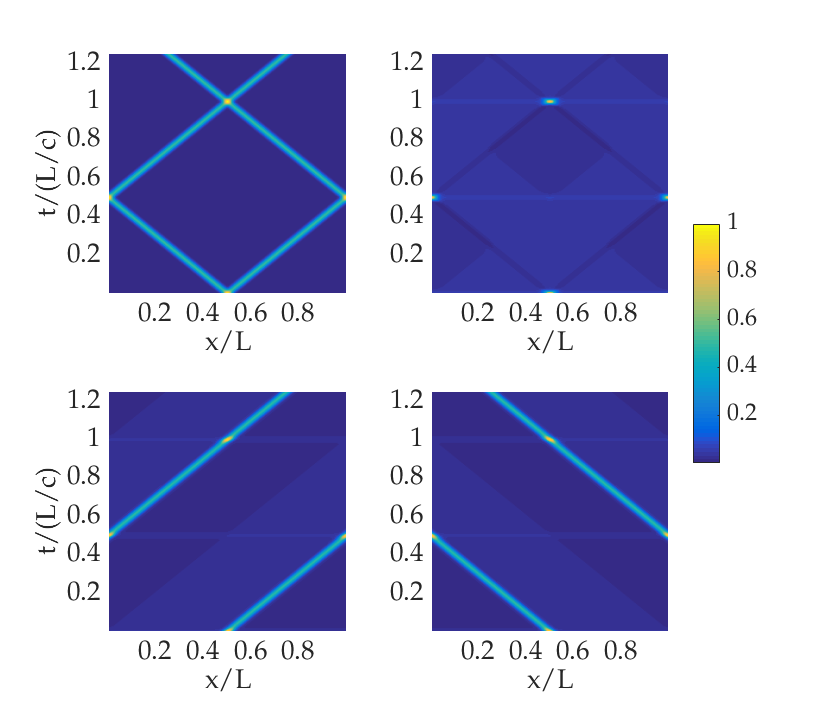}
\caption{
The density of two transported quantities, the left and right going wave: \MSR (top left: Full-order solution, top right: Error of \eqref{eqn:naiveRecon}, bottom left: Right going component of \MSR, bottom right: Left going component of \MSR)
}
\label{fig_naive}
\end{figure}

To illustrate the new approach for problems with multiple transports, we  consider the solution \eqref{solRiemann} of the wave equation \eqref{linearWave} with
\begin{equation*}
q_+ = q_- = \exp( -(x-L/2)^2/(L/50)^2  ),
\end{equation*}
which describes a pressure pulse in a system initially at rest. 
The analytic solution of the density for these two transported quantities is shown in Fig.\ \ref{fig_naive}, top left.  
Due to the solution structure, cf.\ \eqref{solRiemann}, it is known, that the full information is described by just two modes. Ideally, a model reduction approach should detect these modes.
However, while the solution structure is known in this special case, we refrain from using details of it.
 The desired method should find such modes purely data-based. In this way, such a procedure is expected to work also for non-hyperbolic transport.
For the remainder of this subsection, we assume the transport velocities to be known, while in section \ref{sec_determinVel} methods of determining transport velocities based on the snapshot data are discussed.
\medskip

\noindent To identify good ansatz modes, recall that in the case of one transported quantity it is well described by the first modes in a co-moving system, while a decomposition of structures in a different velocity frame leads to the need of many modes to describe the dynamics reasonably well. Consequently, the first mode has an inferior contribution to the dynamics, if the shift velocity does not agree with one of the transport velocities.
Thus, a naive approach to {\it decompose} different velocity components is given by the following procedure.


\hrule
~\\

{\scshape   Multi-Shift \& Reduce (\MSR) }   \\

\textit{Input}:
\begin{itemize}
\item $j\times n$ snapshot matrix X with  $n$ (temporal) snapshots of $j$ grid points
\item  transport velocities  $c_k$, $k=1,\dots,N_s$, $N_s$ number of shift velocities
\end{itemize}

\textit{Procedure}:
\begin{itemize}

\item[{\it 1.}] \textbf{Shift \& Decompose}

Compute the SVDs
\begin{align*}
  U^k\Sigma^k (V^k)^T &=T^{-c_k}\left(X\right)
\end{align*}
with $U^k\in\mathbb{R}^{j,j}$, $\Sigma^k\in\mathbb{R}^{j,n}$, and $V^k\in\mathbb{R}^{n,n}$ for $k=1,\dots,N_s$. Here, $T^{-c_k}$ denotes the discrete shift operator, cf.\ end of section \ref{sec:oneTransport}.

\item[{\it 2.}] \textbf{Truncate} \\
Approximate the matrix $U^k\Sigma^k(V^k)^T$ by neglecting the singular values $\sigma_{r+1}$, $\ldots$, $\sigma_{min\left(k,n\right)}$,  where $r$ is chosen as small as possible, but as large as necessary to obtain a good approximation. This leads to the approximate SVDs
\begin{align*}
    \tilde U^k\tilde \Sigma^k (\tilde V^k)^T    & \approx  U^k\Sigma^k (V^k)^T
\end{align*}
with $\tilde{U}^k \in\mathbb{R}^{j,r}$, $\tilde{\Sigma}^k \in\mathbb{R}^{r,r}$, and $\tilde{V}^k\in\mathbb{R}^{n,r}$ for $k=1,\dots,N_s$.

\end{itemize}
~
\hrule
~\\[1em]
Altogether multi-shift \& reduce (MS\&R)
\begin{equation*}
\MSR: X \to  \tilde \Sigma^k,\tilde  U^k,  \tilde  V^k
\end{equation*}
produces modes containing parts of the field which can be represented well in the respective velocity frame.
A naive approximation of the original field would then be given by
\beqn
X \approx \tilde X = \sum_{k=1} ^{N_s} T^{c_k} \left(   \tilde U^k  \tilde \Sigma^k (\tilde V^k)^T \right).  \label{eqn:naiveRecon}
\eeqn
If the matrices would contain only information of the associated velocity frame this approximation could be exact.
In general a perfect decomposition is not to be expected and instead parts of the field  are  over-represented.

The result of applying this procedure to the pressure pulse example is shown in Fig.\ \ref{fig_naive}, bottom, where we applied \MSR with $r=1$ in each velocity frame.
\noindent \MSR  leads to an approximate identification of the respective velocity component, since this dominates the respective first mode, while the other velocity component gives a relatively small contribution.
The reconstruction is done according to \eqref{eqn:naiveRecon} and the error is shown  in Fig.\ \ref{fig_naive}, top right.

While the main structures are reproduced well, the result is far from perfect.
Especially, in the region where the left and the right going pulse overlap, a large error is produced.
These strongly localized structures with respect to space and time are represented well in both velocity frames, since
a point structure is nearly unchanged by a shift  and can be described perfectly by a dyadic product.
This ambiguity impedes the separation of the two transported quantities.  
Furthermore, the quality of the solution cannot be improved by adding modes, since this counteracts the separation.

To improve the decomposition, a new algorithm is presented in the following which is referred to as the sPOD algorithm.
The \MSR identifies  by construction low rank structures in each velocity frame.
These are used to create a set of ansatz modes to optimize the approximation (\ref{multiAnsatz}) of the field in a least squares sense. 
To this end, the error or residual is decomposed by \MSR to identify the part of the error, which can be represented in each velocity frame.
Finally, the decomposition is  iteratively improved. 
An implementation of the sPOD algorithm and the numerical examples considered in this paper are available from \cite{ReiS18}.

%

~\\
\hrule
~\\

{\scshape The sPOD algorithm}   \\

\textit{Input}:
\begin{itemize}
\item $j\times n$ snapshot matrix X with  $n$ (temporal) snapshots of $j$ grid points
\item  transport velocities  $c_k$, $k=1,\dots,N_s$, $N_s$ number of shift velocities
\end{itemize}

\textit{Procedure}:
\begin{itemize}
\item[{\it 1.}] Initialize
$j\times n$ matrix $ \tilde X = 0$,\;  $j\times r$ matrix $\tilde U^k=0$,\;  $n\times r$ matrix $\tilde V^k  = 0  $
\item[\textbf{loop}]
\item[{\it 2.}]
calculate residual  $$ R = X - \tilde X $$
\item[{\it 3.}] use \MSR\ to create ansatz modes
$$
\MSR:\;  R  \to  \Sigma^k_r, U^k_r,  V^k_r
$$
\item[{\it 4.}]
optimize  in a least squares sense
$$ \min_{\alpha^k, \alpha_r^k  } \left\lVert \hat X - X  \right\rVert_2^2 $$
with
\begin{equation}
\label{eq:sPODansatz}
\hat X = \sum\limits_{k=1}^{N_s}T^{c_k} \left(  \tilde U^k \alpha^k (\tilde V^k)^T  +  U^k_r \alpha^k_r (V^k_r)^T \right)
\end{equation}
where $\alpha^k$ and $\alpha^k_r$ are diagonal coefficient matrices of appropriate dimensions
\item[{\it 5.}]
calculate new modes by SVD
$$    U^k \Sigma^k (V^k)^T =  \tilde U^k \alpha^k (\tilde V^k)^T  +  U^k_r \alpha^k_r (V^k_r)^T $$
and truncate it as in step 2 of \MSR
$$  \tilde U^k \tilde \Sigma^k (\tilde V^k)^T  \approx U^k \Sigma^k (V^k)^T     $$
\item[{\it 6.}]
update approximation
$$ \tilde X = \sum\limits_{k=1}^{N_s}T^{c_k} \left(\tilde U^k \tilde \Sigma^k (\tilde V^k)^T\right) $$
\item[\textbf{until}] $\left\lVert R\right\rVert_2$ does not reduce further
\vspace{0.1cm}

\end{itemize}

\hrule
~\\[2em]

Due to the initialization, the first decomposition is \MSR of the original field.
In further steps the residual is reduced by providing modes which are constructed by \MSR of the residual and allow to remove structure which is multiply accounted for.
For example the strong peaks in the residual in Fig.\ \ref{fig_naive} are prominently visible in the residual modes $U^k_r$ and can be used to remove this error to a large degree.
The iteration in sPOD offers the possibility to remove it up to any achievable tolerance.

An open question is how to choose the $r$ in the \MSR for each frame in an optimal way.
For cases where $r$ is not clear from physical considerations, we propose two different heuristics. 
\label{modeNumHeu} 

One method is to perform the sPOD decomposition with a large number of modes per frame and select the most important ones afterwards. 
This is done by sorting the singular values of all reference frames in one list and choosing the modes associated with the largest singular values. 
If a certain approximation error is prescribed, one has to calculate the residual from the reconstruction.  
It cannot, however, be calculated from  the singular values directly, due to the  non-orthogonality of modes belonging to different frames.

An alternative way of choosing $r$ is to start with small numbers of modes and successively add modes in a greedy fashion.
We explain this in more detail with an example of three reference frames and starting with zero modes for each frame, i.\,e., $\boldsymbol{r}_0=[0\,0\,0]$, where $\boldsymbol{r}_k$ denotes the vector containing the numbers of modes of each velocity frame after the $k$th greedy iteration.
At the first iteration, different sPOD approximations are computed for $\boldsymbol{r}=[1\,0\,0]$, $\boldsymbol{r}=[0\,1\,0]$, and $\boldsymbol{r}=[0\,0\,1]$.
The errors of the three different sPOD approximations are compared and we add only one mode to that frame which corresponds to the smallest error, for instance $\boldsymbol{r}_1=[0\,1\,0]$ if the comparison shows the greatest advantage of adding a mode to the second frame.
This procedure can be continued until a certain error tolerance is achieved and can be implemented as a loop around the sPOD algorithm.
Of course, the greedy algorithm is based on a locally optimal choice which does not guarantee that the resulting vector $\boldsymbol{r}$ is optimal.
Furthermore, this approach can be quite costly if the amount of data and the number of frames is big since the sPOD algorithm has to be performed multiple times. 
Nevertheless, in the considered numerical examples the greedy algorithm yields sPOD approximations with high accuracy while only using a small number of modes per frame.

\begin{figure}
\centering
\includegraphics[scale=0.23]{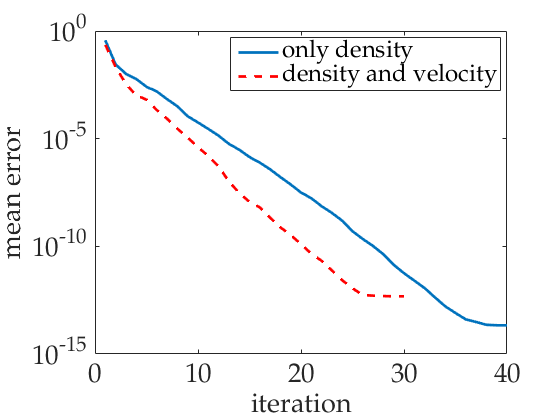}
\includegraphics[scale = 0.23,clip=true,trim = 2.5cm 0cm 1.5cm 0cm]{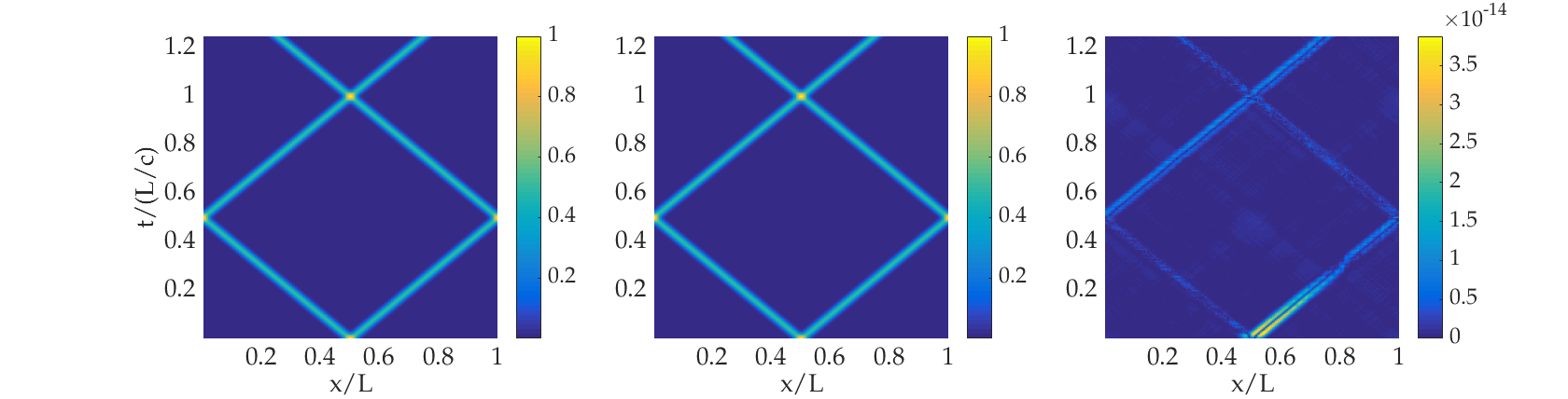}
\caption{Two transported quantities. Left: Convergence of the sPOD algorithm (2 modes). The error of the sPOD approximation decreases to less than $5\times10^{-13}$. Right: Full-order solution of the density, sPOD approximation, and the error (note the color scale).}
\label{fig_sPOD1D}
\end{figure}

The sPOD algorithm is tested for the density of the acoustic pulse with one mode per frame.
 The convergence behavior can be seen from the solid graph of Fig.\ \ref{fig_sPOD1D}, left, 
 where the convergence is indicated by the decreasing mean error over iterations.
The considered error measure is given by the $2$-norm of the error divided by the $2$-norm of the full-order snapshot matrix, i.\,e.,
\begin{equation}
\label{eq:relMeanError}
\mathrm{mean\;error} = \frac{\left\Vert X-\tilde{X}\right\Vert_2}{\left\Vert X\right\Vert_2}.
\end{equation}
Some remarks follow.

\paragraph{1} While in the first iterations the sharp structures allow a good separation of the different transport directions, the residual tends to lose this structure.
The interference of the components of different velocity create a broad structure, visible as the broad diagonal stripes in Fig.\ \ref{fig_naive}.
This tends to slow the convergence after the first iterations.

\paragraph{2}\label{linearchange} For the considered example, the analytic solution has the form of two modes, each in a different frame, i.\,e., $X=T^{c_1}(u^1\sigma^1(v^1)^T)+T^{c_2}(u^2\sigma^2(v^2)^T)$. This allows to consider
the modes from the initial \MSR as the modes of the exact solution perturbed by the component of the respective other velocity frame. For example for the first frame of reference, the shifted snapshot matrix is
\begin{equation*}
T^{-c_1}(X)=u^1\sigma^1(v^1)^T+T^{c_2-c_1}(u^2\sigma^2(v^2)^T).
\end{equation*}
Assuming the perturbation $T^{c_2-c_1}(u^2\sigma^2(v^2)^T)$ to be small, the perturbed singular values and vectors can be expressed as a linear perturbation of the unperturbed ones with $\delta \sigma^1 = (u^1)^T \delta X^1 v^1 $ and $  \delta u^1 \sigma^1  = \delta X^1v^1 - u^1 \delta \sigma^1 $ characterizing the linear change of $\sigma^1$ and $u^1$.
For the considered pressure pulse example, the mode obtained by \MSR is well approximated by $\tilde u^1\approx u^1+ \delta u^1$. This suggests that the exact mode is indeed dominant in its reference frame and the component of the other frame is approximately a linear perturbation.
The broad stripes in Fig.\ \ref{fig_naive} are thereby explained  as  an overlap of the exact mode and the mode of the other frame in the respective velocity frame.

\paragraph{3} Broad structures in the residual can be better resolved with more modes. Therefore, increasing the number of modes $r$ in the \MSR for the residual tends to result in an improved convergence. Indeed, we have observed for the wave equation with some less sharp initial conditions that choosing $r=1$ prevents the procedure from converging to the analytic solution. For these examples increasing the number of residual modes $r$  improves convergence and results in finding the analytic solution up to any given achievable tolerance. However, there is a trade-off between convergence and computational effort both of which have to be accounted for when choosing the number of residual modes.

\paragraph{4} The ansatz in \eqref{eq:sPODansatz} can also be extended so as to include cross terms of the form $U^i(V^j)^T$ with $i\neq j$.
This can be achieved by optimizing over dense coefficient matrices $\alpha^k$ and $\alpha^k_r$ rather than over diagonal ones.
We have observed that this often leads to a faster convergence of the sPOD iteration which is due to the increased number of degrees of freedom for the optimization in step 4 of the sPOD algorithm.
Nevertheless, we omit the cross terms for reasons of computational effort.
The number of unknowns in the resulting system of linear equations to be solved in step 4 would scale with $N_s r^2$ instead of $N_s r$.
This is a minor issue for examples in one dimension but it leads to high computing times when dealing with the two-dimensional example presented in section \ref{sec_twoD}.

\paragraph{5} As an alternative to reducing the different velocity frames simultaneously, one could proceed sequentially by removing the detected structures in a greedy fashion.
This is done in the recent work \cite{RimMoeLeVeque2017}. It is clear from Fig.\ \ref{fig_naive} that after the first iteration each of the modes contains structure from the respective other velocity frame. However, due to the greedy character these artifacts are not removed, since the modes from the first iteration remain unchanged. As a consequence, the sequential procedure does not reduce to descriptions with just one mode per velocity frame in the pressure pulse example.
Moreover, the obtained modes do not reflect the physics of two moving waves properly, which might also be disadvantageous for  model reduction.
\smallskip
%

The error of approximation obtained by the sPOD algorithm after 40 iterations is shown in Fig.\ \ref{fig_sPOD1D}, right. As can be seen, the solution is approximated excellently by just two modes. This is the desired low-dimensional representation we have been looking for, agreeing with the  analytic solution given by two modes. Note that, since only the density was provided, the Riemann invariants cannot be calculated, suggesting that the construction did not implicitly use the hyperbolic structure of the equation.


\begin{figure}
\centering
\includegraphics [scale=0.3]{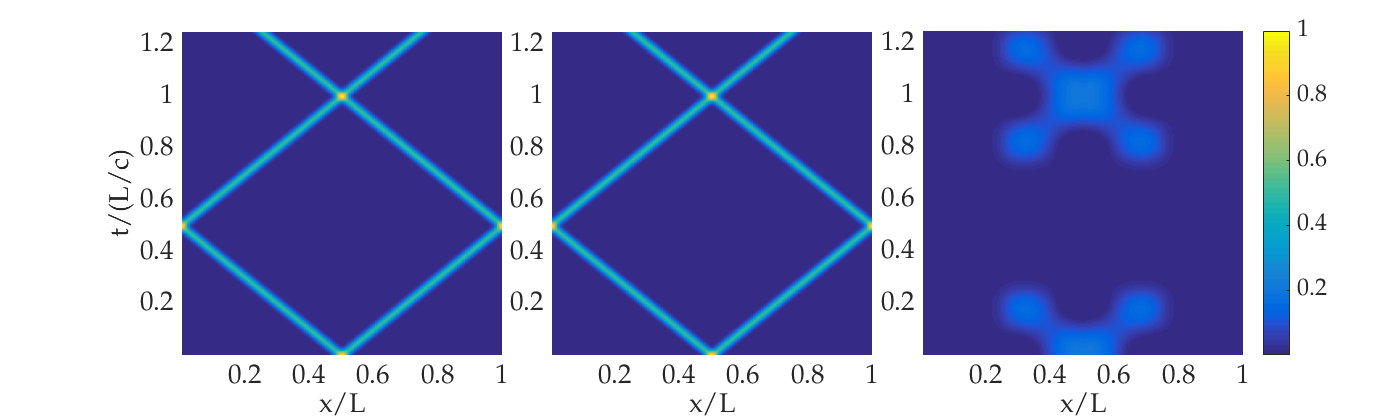}
\caption{Two transported quantities (only density): Comparison of full-order solution with sPOD and POD approximation (left: Full-order solution, middle: sPOD approximation (2 modes), right: POD approximation (2 modes))}
\label{fig:1D_density_POD_vs_sPOD}
\end{figure}

The good performance of the sPOD algorithm gets more striking when comparing it to the standard POD.
For this purpose, Fig.\ \ref{fig:1D_density_POD_vs_sPOD} depicts the comparison between the full-order solution, the sPOD approximation with two modes, and the POD approximation with two modes. It is obvious that the POD approximation is highly inadequate, whereas the sPOD approximation matches the full solution excellently.
To obtain the same accuracy as the two sPOD modes (error less than $3\times 10^{-14}$, cf.\ blue, solid graph in Fig.\ \ref{fig_sPOD1D}, left), more than $80$ POD modes are required. Furthermore, if only two POD modes are used, the relative mean error, as defined in \eqref{eq:relMeanError}, is almost $1$. Of course, this insufficient performance of the POD was to be expected, since deliberately we have chosen an example, which provides a big challenge for the POD. Nevertheless, this test case which is generic for many practical problems gives a first impression of the potential of the proposed sPOD.

%

For all the considerations made so far, the reduction was based on the density alone. When considering the density and the velocity together, the convergence of the sPOD algorithm appears to be faster even though the eventually achieved error is somewhat higher but the order of magnitude is $10^{-13}$ (see Fig.\ \ref{fig_sPOD1D}, left). 
Apart from the convergence behavior, the results for the case that density and velocity are considered together are very similar to the results shown before and omitted therefore.

\subsection{Determination of the shift velocities}
\label{sec_determinVel}
If the  shift velocities are unknown, different strategies can be used to obtain good candidates for them.
This is largely independent of the sPOD, which is the main focus of this paper, but important for practical use.
Often physical insight reveals  (some of) the involved  velocities, i.\,e., the flow or sound velocity.
The main target of the sPOD is the decomposition of strongly transport-dominated cases. This class of problems is especially interesting since classical methods like the POD fail here most severe, see e.\,g{.} \cite{CagMS16}. In these cases the
transport velocities can often simply be determined by usage of data-based tracking of peaks or threshold values within the snapshot matrix.
Another method is discussed in the following where the velocity detection is performed via a maximization of singular values.
The tracking-based method is applied in the example of section \ref{sec:twoShocks} while the velocities of the example in section \ref{sec_twoD} are determined by means of the singular value maximization.

\begin{figure}
\centering
\includegraphics [scale = 0.27,clip=true,trim = 0cm 0cm 0cm 1cm]{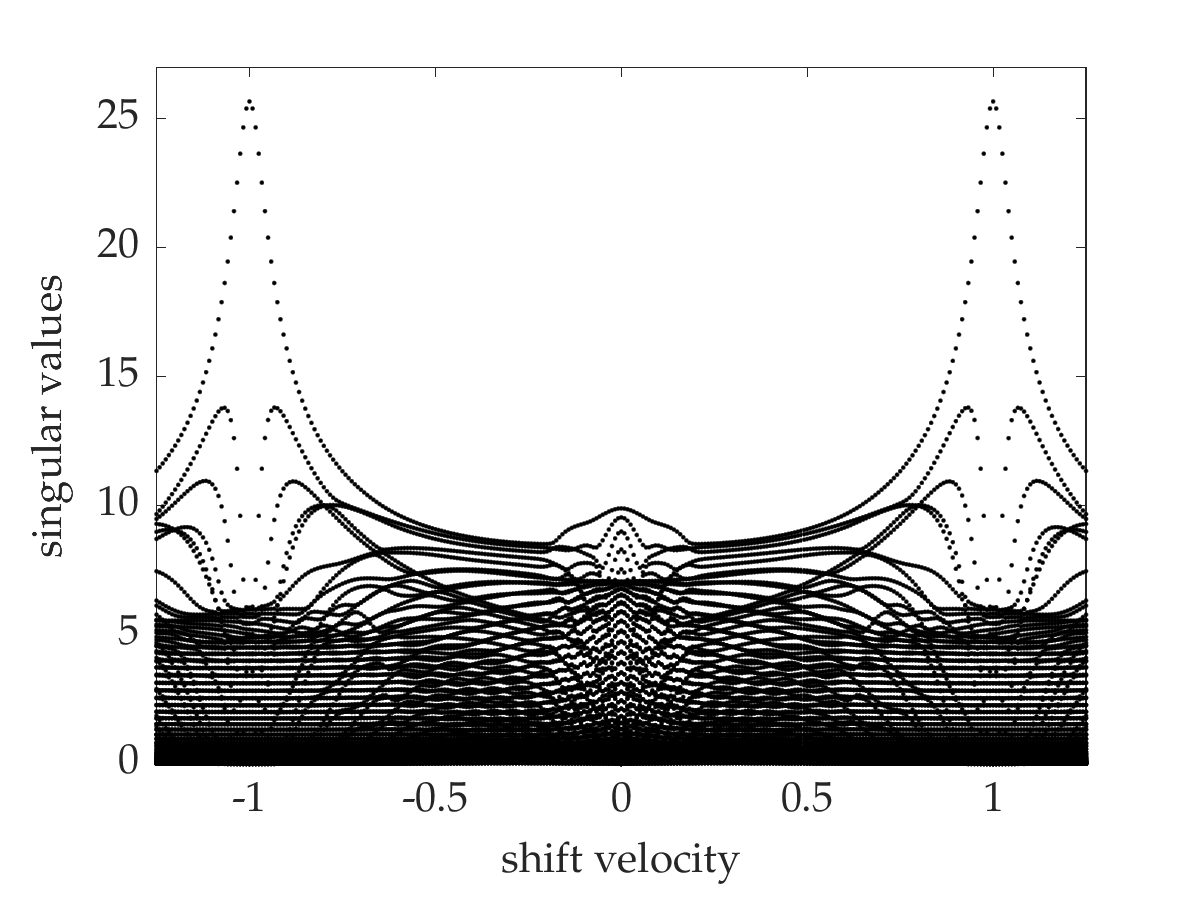}
\includegraphics [scale = 0.28,clip=true,trim = 0cm 1cm 0cm 0cm]{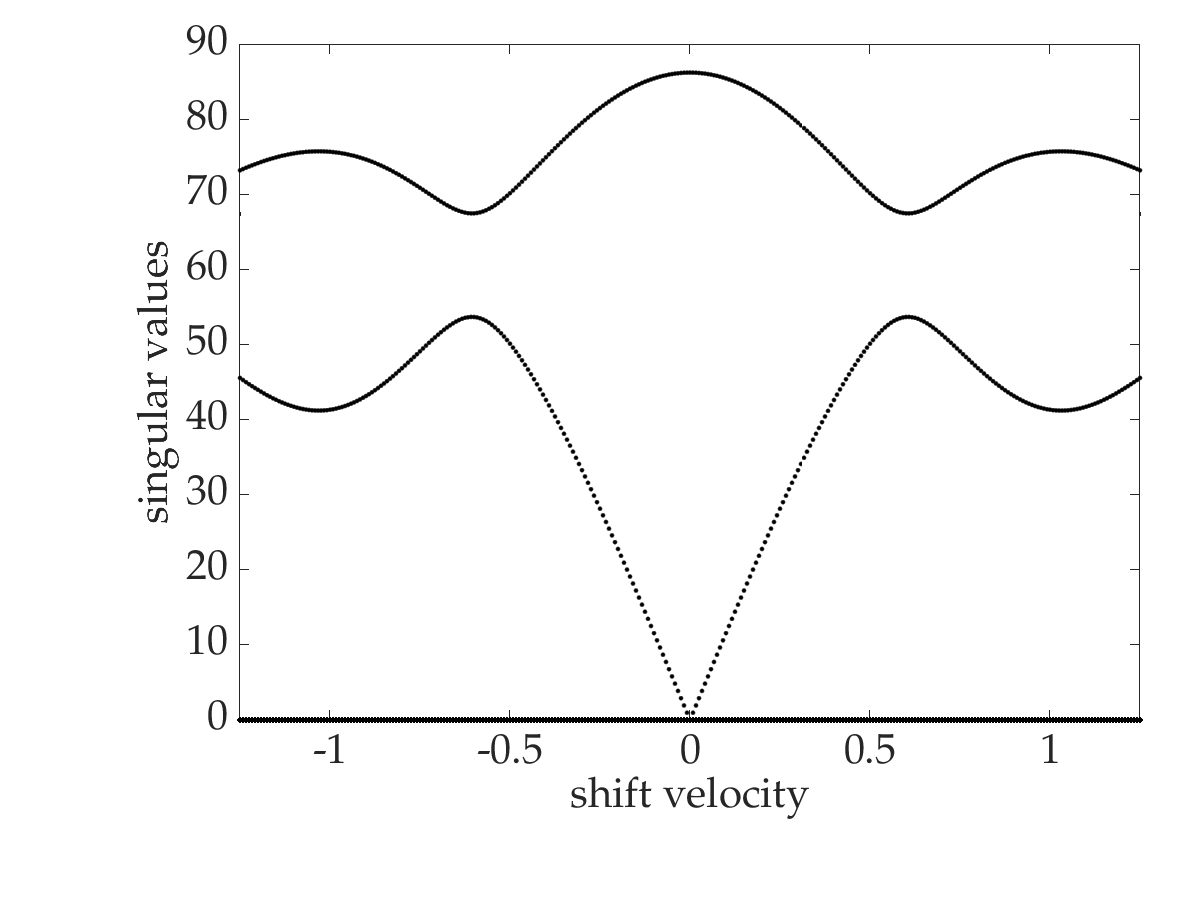}
\caption{
The singular value spectrum of the density snapshot matrix as a function of the shift velocity.
Left:  The  pressure pulse, which creates two traveling waves with the speeds $c_\pm = \pm 1$, which are clearly visible as maxima of the leading singular value.
Right: The standing wave. The zero velocity has the biggest leading singular value. The positive and the negative sound speed are also visible as local maxima.
}
\label{fig_svd_oneD_trans}
\end{figure}

A purely data-driven method  is obtained by  examining the singular value spectrum of the shifted snapshot matrix as a function of the shift velocity, i.\,e.,
 the (constant) velocity the time-dependent shift is based on. 
In Fig.\ \ref{fig_svd_oneD_trans}, left, the spectrum of shifted snapshot matrices of the example considered in section \ref{sec:multipleTransports} is shown for a range of shift velocities between $-1.25$ and $1.25$.
First, it should be noted that, assuming periodic boundaries, the square integral of the solution
\begin{equation*}
\int\limits_0^{t_{end}} \int\limits_0^L  (q(x,t))^2  \mathrm{d}x \mathrm{d}t
\end{equation*}
does not change by a shift in $x$-direction.
A numerical approximation of the shift  keeps this invariance up to the interpolation error.
Consequently, the Frobenius norm, which is directly connected with
the singular values by
\begin{equation*} \left\Vert X\right\Vert_F^2 = \sum_{i,j} X_{i,j}^2 = \sum_i \sigma_i^2,
\end{equation*}
is also shift-invariant up to the interpolation error.
This allows to directly compare spectra for different shift velocities. If a certain singular value increases by a change of shift velocity, others have to decrease to keep the sum of the squared singular values constant. The transport velocities, which are the positive and the negative sound velocity $c_\pm=\pm 1$, are clearly visible from the maxima of the leading singular value.
The classical POD is recovered for a shift velocity of zero, where a slow decay of the singular values can be seen in Fig.\ \ref{fig_svd_oneD_trans}, left.
The shift dependence of the pressure pulse case is  contrasted with the case of a standing wave, Fig. \ref{fig_svd_oneD_trans}, right.
Here the zero shift velocity leads to a maximization of the leading singular value.
However, the sound velocities are also visible as local maxima.
The standing wave can be represented by two traveling waves with velocities $\pm c$, so that these are also reasonable candidates for shift velocities.
It should be noted that the spectra are only a first indicator.
If multiple transport velocities are involved, each of the transported quantities influences the singular spectra of the others.
For instance the pressure pulse can be expressed by one mode per transport velocity, which is not obvious from Fig.\ \ref{fig_svd_oneD_trans}, left,
since each of the transported quantities slows the decay of the singular spectrum of the other one.
The simple sampling used here can be replaced by a gradient method, since a change of singular values with respect to a change of the matrix can easily be calculated, similar to the considerations on page \pageref{linearchange}, remark 2.

%

\subsection{Uniqueness and robustness of the method}

For the case $N_s=1$ and $c_1=0$ in \eqref{eqn:multi_struc}, the sPOD reduces to the classical POD and hence for this special case an optimal approximation is given due to the optimality of the SVD.
A simple counterexample shows that such a global optimality cannot be guaranteed for the sPOD in general.
Consider again the pressure pulse example from section \ref{sec:multipleTransports} and add a third velocity $c_3=0$ to the sPOD algorithm.
Since the solution can be described perfectly by the two moving modes, the zero velocity frame is redundant.
However, we can describe the analytic solution also by adding a harmonic mode to the zero velocity frame and compensate for this by adding appropriate components to the two co-moving frames, cf.\ \eqref{waveRI} and \eqref{solString}. 
Accordingly, the analytic solution allows different, but equivalent low-dimensional representations.  
Adding a mode in the zero frame and subtracting from the Riemann invariants does not change the composed solution, e.\,g., 
	
%
\begin{align*}
q(x_i,t_n) &= T^{c}\left(\begin{bmatrix}
\rho_0\\
c
\end{bmatrix}q_{+}\left(x\right) \right)+T^{-c}\left(\begin{bmatrix}
\rho_0\\
-c
\end{bmatrix}q_{-}\left(x\right) \right)\\
 &= T^{c}\left(\begin{bmatrix}
\rho_0\\
c
\end{bmatrix}\left(q_{+}\left(x\right)-\cos\left(k_1x\right)\right) \right)+T^{-c}\left(\begin{bmatrix}
\rho_0\\
-c
\end{bmatrix}\left(q_{-}\left(x\right)-\cos\left(k_1x\right)\right) \right)\\
&+2T^{0}\left(\begin{bmatrix}
\rho_0\cos\left(k_1ct\right)\cos\left(k_1x\right)\\
c\sin\left(k_1ct\right)\sin\left(k_1x\right)
\end{bmatrix} \right).
\end{align*}
 This shows that the separation is not unique and all of these solutions are fix points of the sPOD iteration.
Indeed, the numerical experiments show that the sPOD algorithm does not remove the component of the zero velocity frame but instead converges to a solution which is represented as a sum of the three provided frames. The accuracy is still comparable to the sPOD just using the two co-moving subspaces.
This means that the sPOD is not optimal in the sense that a minimum total number of modes is not guaranteed. Note that a decomposition based on three modes is still much better than the classical POD in this example. One could think of strategies to enforce removing redundant information, for instance by augmenting the least squares approach in step 4 of the sPOD algorithm with an $\ell_1$-norm regularization (see e.\,g{.} \cite{BoyV09}). However, this exceeds the scope of this paper.

If the transport velocity is not determined correctly but with a small error, an increased number of modes is required to obtain the same accuracy as with the correct transport velocity. However, as can be seen in Fig.\  \ref{fig_svd_oneD_trans}, the leading singular value decreases continuously from the maxima.
Hence, small errors from the velocity detection lead to a small increase of the required number of modes.
This growth depends essentially on the width of the maximum.
A sharp maximum, i.\,e., a maximum exhibiting a strong decay of the leading singular value close to it, requires a precise velocity determination. On the other hand, sharper maxima are also easier to determine, which suggests that also for these cases a sufficiently accurate determination of the dominant transport velocity is possible with the presented method.
\subsection{Two crossing shocks}\label{sec:twoShocks}

In this section we consider a more physical and challenging example: The crossing of two shock waves.
The example of two crossing shocks is a generic phenomenon which occurs in many applications. 
The non-constant shock velocities and amplitudes are inherently nonlinear and change during the crossing of the shocks.
However, this is not a problem for the sPOD decomposition which can also work with variable velocities.
To this end, the shift operator introduced in section \ref{sec:oneTransport} is generalized to account for time-dependent velocities by replacing the linear expression $ct$ by a general time-dependent shift coordinate $x^{sh}\left(t\right)$, i.\,e., $\mathcal{T}^{x^{sh}} f\left(x,t\right) \vcentcolon = f\left(x-x^{sh}\left(t\right),t\right)$.

The Euler equations for an ideal gas with constant heat capacity are simulated by the skew-symmetric scheme with a shock filter, as  described in \cite{ReiS14}.
The nonlinear shock filter  \cite{BogeyCacquerayBailly2009} chooses the filter strength on the base of a shock detector.
The initial condition is $\rho = [1.0138,0.6000 , 1.0138 ]$, $u=[ -96.2197,-310.0174, -523.8264]$, $p = [ 1.2720, 0.6000, 0.2720] $, yielding two shocks colliding with shock Mach-number of $Ma=1.4$.
The shocks are at $x= [1/3, 2/3]$. The initial shock structure is chosen to have a near steady state with the given shock filter.
The adiabatic exponent is $\gamma=1.4$.

\begin{figure}
\centering
\includegraphics[scale = 0.3,clip=true,trim = 2.5cm 0cm 1.5cm 0cm]{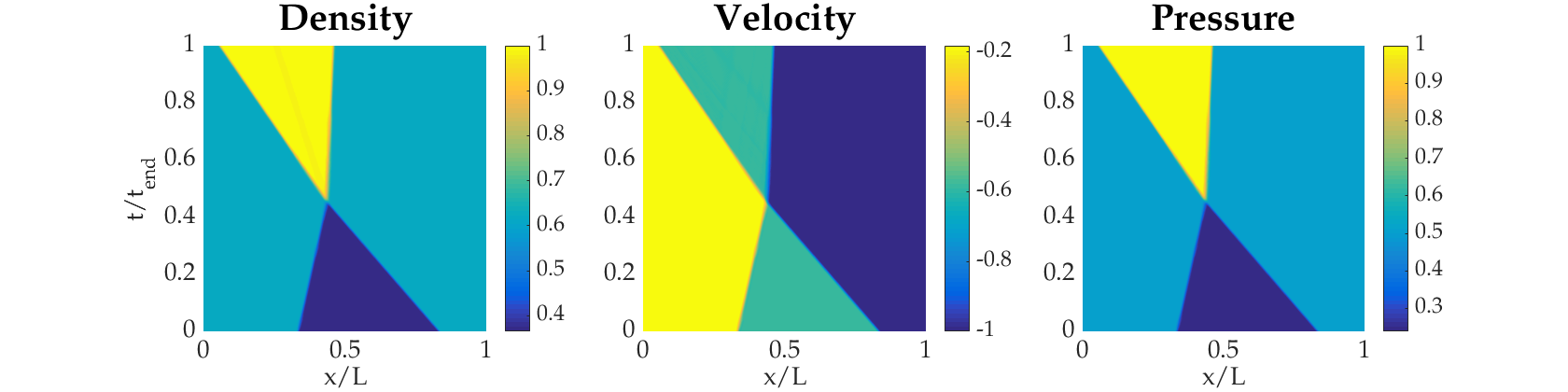}
\caption{Full-order solution of two crossing shocks: Pseudocolor plots of the density (left), velocity (middle), and pressure (right)}
\label{fig:twoShocks_FOM}
\end{figure}

The simulated solution is depicted in Fig.\ \ref{fig:twoShocks_FOM} for the density, the velocity, and the pressure.
The depicted space-time diagrams represent the transposes of the respective snapshot matrices which are used as an input of the sPOD algorithm.
The snapshots have been normalized such that the highest occurring absolute value of each quantity is equal to one to avoid numerical errors due to the different scales. 
In contrast to the examples before, the case of two crossing shocks exhibits non-periodic boundaries. This leads to the situation that we need values from outside of the computational domain when shifting the snapshot matrix. We treat this here by a constant extrapolation over the boundaries, while a general treatment of non-periodic boundaries within the sPOD framework will be addressed in a future work. 

We are looking for a decomposition of the snapshot matrices with a prescribed relative error tolerance of one percent.
To apply the sPOD the time-dependent shift needs to be  determined.
Here, we apply a simple threshold search within the snapshot matrix of the velocity.
More precisely, we define two threshold values characterizing the velocity jumps at the two shocks.
To this end, the respective threshold value has to be chosen in between the constant values at the left- and on the right-hand side of the respective shock.
Regarding the velocity snapshots depicted in Fig.\ \ref{fig:twoShocks_FOM}, we have chosen the thresholds $-390m/s$ (normalized $-0.74$) for the right border of the zone between the shocks (border between cyan and dark blue in middle plot of Fig.\ \ref{fig:twoShocks_FOM}) and $-200m/s$ (normalized $-0.38$) for the left border (between yellow and cyan).
Based on these thresholds, the time-dependent shifts are determined by searching the last value, which is larger than the corresponding threshold, in each column of the snapshot matrix (each discrete time step).

\begin{figure}

\centering
\begin{subfigure}[b]{0.49\textwidth}
\centering
\includegraphics[scale=0.42,clip=true,trim=0cm 0cm 1cm 0cm]{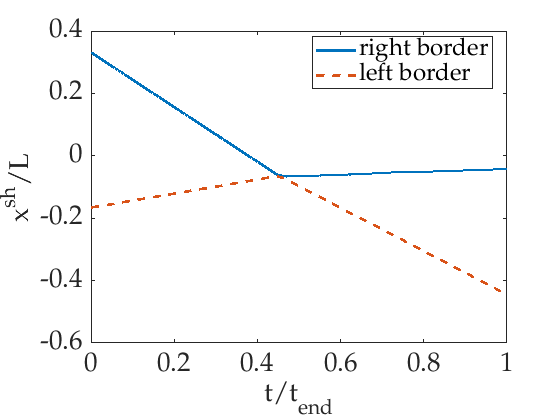}
\caption{Shift coordinates determined via threshold search}
\label{fig:twoShocks_cs1}
\end{subfigure}
\begin{subfigure}[b]{0.49\textwidth}
\centering
\includegraphics[scale=0.42,clip=true,trim=0cm 0cm 1cm 0cm]{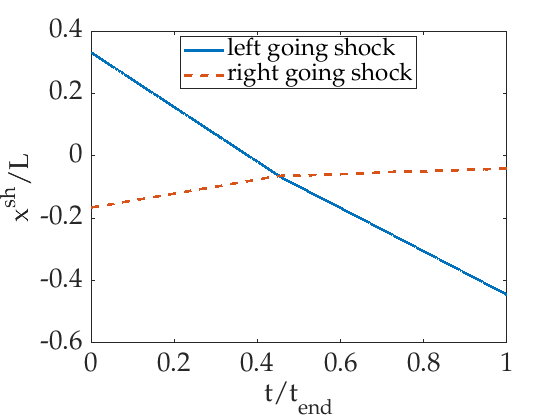}
\caption{Shift coordinates after correction to account for the crossing}
\label{fig:twoShocks_cs2}
\end{subfigure}
\caption{Shift coordinates}
\label{fig:twoShocks_cs}
\end{figure}

The shift coordinates obtained from the threshold search do not intersect but rather bounce off each other, see Fig.\ \ref{fig:twoShocks_cs1}.
In order to account for the crossing of the shock waves, the parts after the time of shortest distance between the two curves are switched and in the intermediate area a linear interpolation is applied.
The resulting shift coordinates, cf{.} Fig.\ \ref{fig:twoShocks_cs2}, provide the basis for the shift matrices applied in the sPOD algorithm.
In order to satisfy the given accuracy requirement of one percent, we apply the sPOD algorithm together with a greedy increment of the number of modes as explained in section \ref{sec:multipleTransports}.

\begin{figure}
\centering
\includegraphics[scale = 0.3,clip=true,trim = 2.5cm 0cm 1.5cm 0cm]{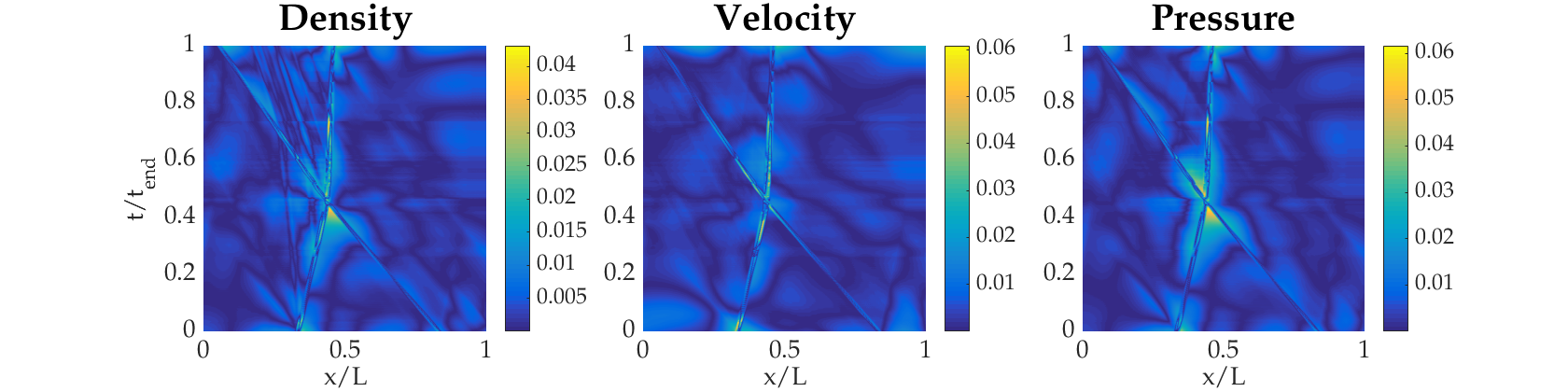}
\caption{Error of the sPOD approximation of two crossing shocks: Pseudocolor plots of the density (left), velocity (middle), and pressure (right)}
\label{fig:twoShocks_error}
\end{figure}

We end up with seven modes, four for the left going and three for the right going wave.
The assembled sPOD approximation is nearly indistinguishable from the simulated solution,  
the error plots are depicted in Fig.\ \ref{fig:twoShocks_error}.
The approximation agrees neatly with the full-order solution which can be comprehended by looking at the error whose maximum amplitude is around six percent of the maximum amplitude of the full-order solution.
The relative mean error is less than one percent as required.
Note, that this could be achieved with just seven modes in total whereas a POD approximation needs 51 modes to attain the same accuracy. 

The error plotted in Fig.\ \ref{fig:twoShocks_error} reveals structures. 
First, the shocks are visible.
This is caused by a not strictly constant shock structure due to small variations in the strength of the shock adaptive filter.
The dynamics of this numerical artifact is not fully described by the used number of sPOD modes.
Second, waves in the density are emerging from the shock crossing and from the initial position of the left shock. 
The two extra waves travel with the flow velocity.
The steady state condition of the left shock is slightly perturbed.
Small perturbations travel as characteristic waves, so that this perturbation creates a so-called entropy wave.
A similar perturbation is created by a change of shock filter strength during the shock crossing.
 These structures in the error could be removed by adding frames with the flow speeds to the decomposition. 

\section{A two-dimensional model problem\label{sec_twoD}}
\label{sec_2D}

In the last section we have considered one-dimensional test examples, which provide a big challenge for the classical POD.
The sPOD performs excellently, because its structure is adequate for the transport-dominated examples.

\begin{figure}
\centering
\includegraphics[scale = 0.4,clip=true,trim = 2.5cm 0cm 2cm 0cm]{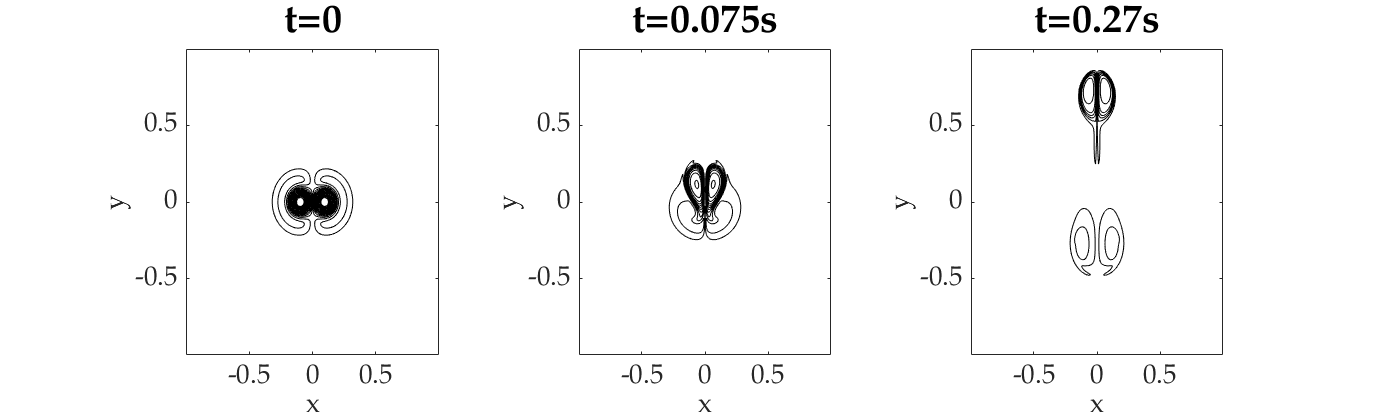}
\caption{2D vortex pair: Full-order solution (vorticity contour lines at $-270$, $-250$, $\ldots$, $270$)}
\label{fig:vortexPair_FOM_Solution}
\end{figure}

In this section we want to explore the capabilities of the introduced sPOD algorithm by considering a two-dimensional, transport-dominated, non-hyperbolic example with non-trivial velocities.
We consider a flow governed by the incompressible Navier-Stokes equations with a vortex pair as  initial condition.
The initial conditions are created from two single vortices with vorticity
\begin{equation*}
\omega_{0,i} = \omega_{e,i} \left(1- (r_i/r_0)^2\right)\,\exp\left(-(r_i/r_0)^2\right)
\end{equation*}
where $i=1,2$, $r_i = \sqrt{(x-x_{0,i})^2+(y-y_{0,i})^2   } $ is the distance from the respective vortex core and $r_0$ denotes the vortex size.
The size of each vortex is chosen to be $r_0 = 0.1$ and the centers are at $(x_{0,1/2},y_{0,1/2}) = (\pm 0.1 , 0)$.
The strengths are $\omega_{e,1/2} = \mp 299.5$ and the viscosity of the fluid is $\nu = 1/Re= 1/1000$.
We consider a  periodic domain to avoid the further complication of the
  boundary treatment for the sPOD. 
 The dynamical behavior of the vortex pair is simulated by means of the energy conserving, skew-symmetric scheme described in \cite{Reiss2015}.
 The periodic domain is discretized by 512$^2$ equidistant points with a time step of $\Delta t = 8\cdot 10^{-4}$s with a fractional-step time stepper.
 The solution is depicted in Fig.\ \ref{fig:vortexPair_FOM_Solution} by means of a contour plot for different times.
The initial vorticity field induces a movement of the vortex pair in positive $y$-direction.
Additionally, a secondary, weaker vortex pair moves in negative $y$-direction with a smaller transport velocity.
Due to this secondary vortex, two different, non-trivial velocities are present.
The simulated time is chosen such that the primary and the secondary vortex pair do not meet again in the periodic domain.
Compared with the linear wave equation from section \ref{sec_oneD}, this example is significantly more challenging for several reasons.
First, we do not know the analytic solution of this test case.
Second, the transport velocities are non-constant and unknown a priori.
Third, the two-dimensional problem not only leads to larger data sets, but tests the applicability in more dimensions.

\begin{figure}
\centering
\includegraphics[scale=0.3]{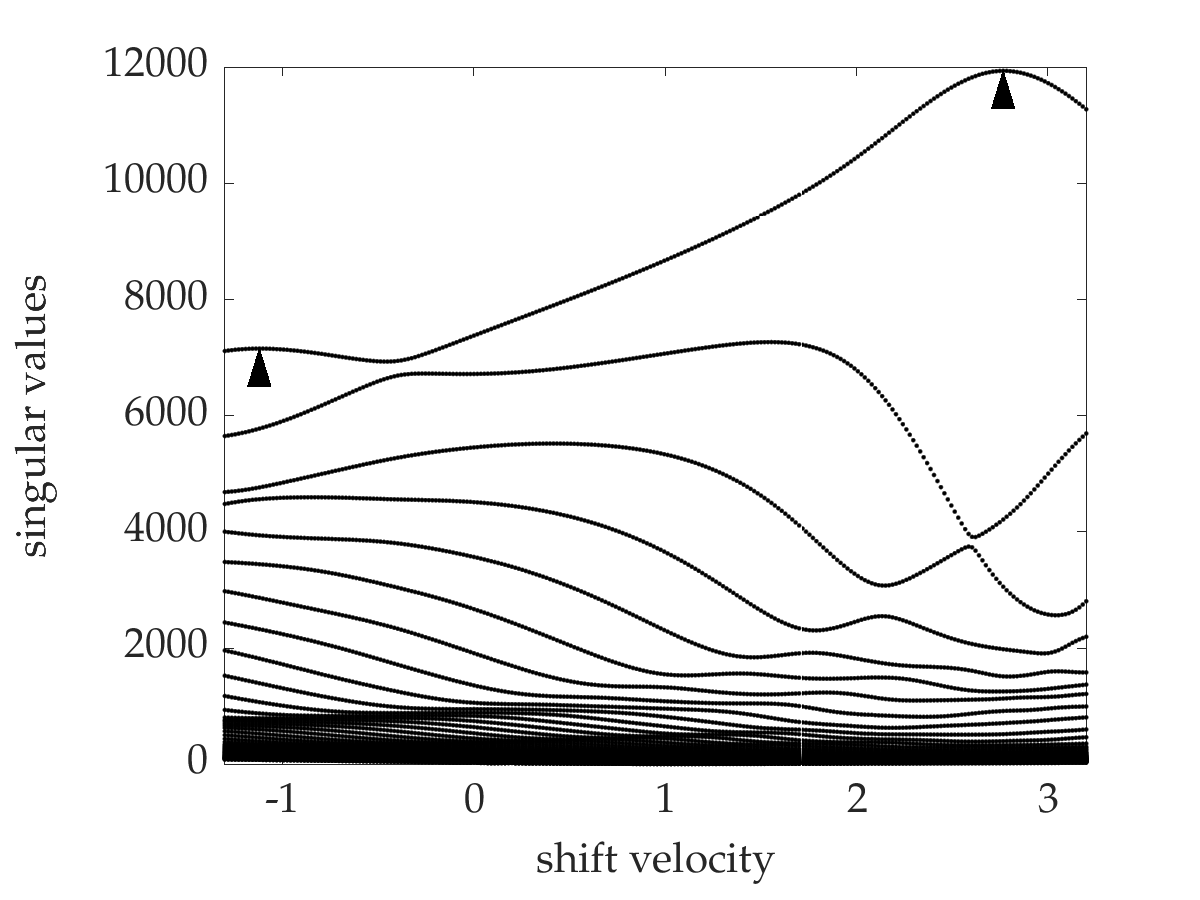}
\includegraphics[scale=0.3]{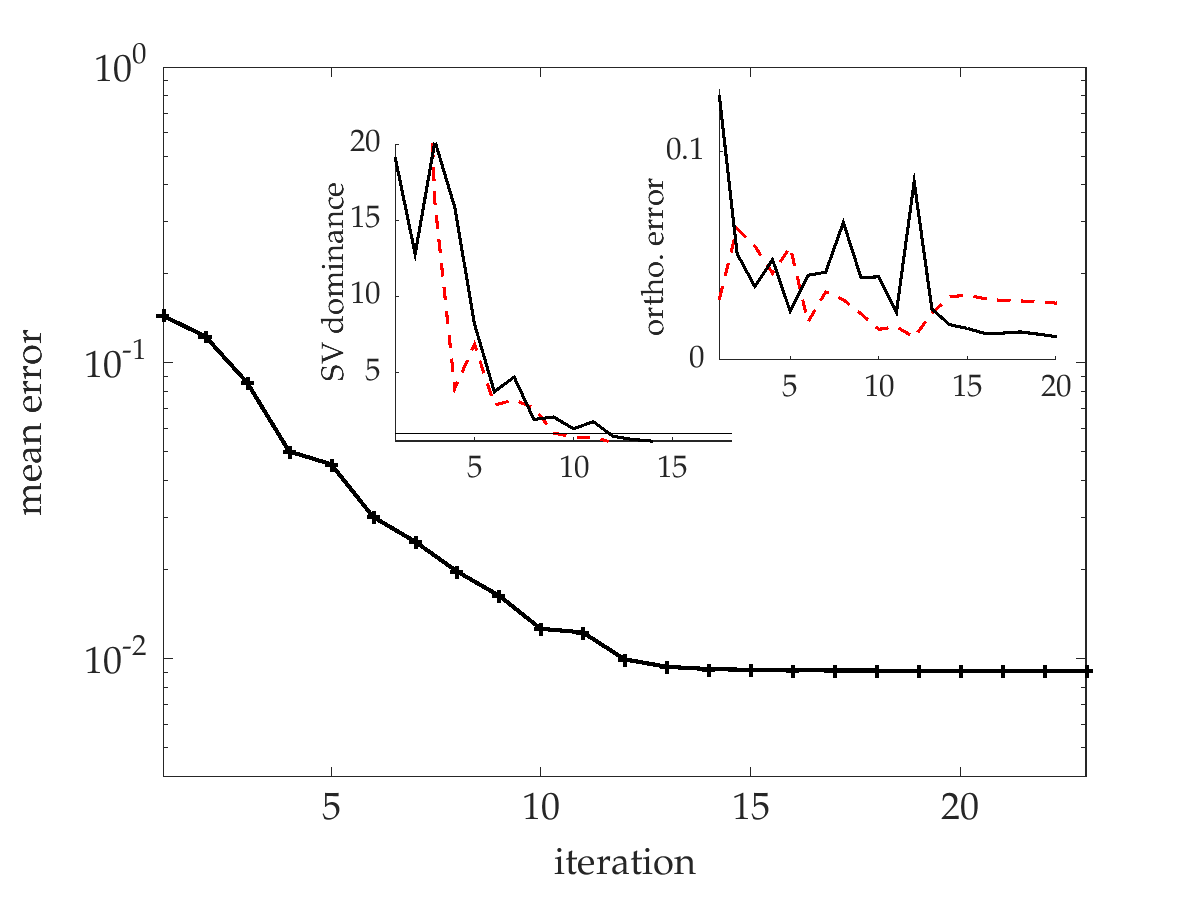}
\caption{2D vortex pair: Left: Velocity scan in y-direction. The two marked maxima at $-1.1193$ and  $2.7663$ are used for the decomposition.
Right: Convergence of the sPOD measured by the norm of the residual related to the norm of the field. In the inset the ratio of the largest singular value of the residual and the smallest singular value of the sPOD approximation in the respective frame is plotted (solid line for the primary, red dashed for the secondary vortex pair).
The horizontal line marks the value one below which the resolved singular values are greater than those of the residual.
The second inset shows the {\it maximal} orthogonality error, cf. \eqref{defConvergence}.}\label{fig:vortexPair_scan}
\end{figure}

Before we apply the sPOD algorithm as introduced in section \ref{sec_oneD}, we need to find proper candidates for the dominant transport velocities. 
We are looking for two velocities which we treat as constant, even though this is not strictly true.
For the velocity determination, we again consider the singular values of shifted snapshot matrices where we only shift in $y$-direction.
Since the considered case is two-dimensional, we need to reshape the snapshot matrix as it would be done for the POD.
The singular values as a function of the transport velocity are shown in Fig.\ \ref{fig:vortexPair_scan}, left.
The maxima are at $-1.1193$ and $2.7663$ which are used as the  transport velocities in the decomposition of the secondary and primary vortex pair, respectively.

\begin{figure}
\centering
\includegraphics[width=0.32\linewidth,clip=true, viewport = 20 0 560 420]{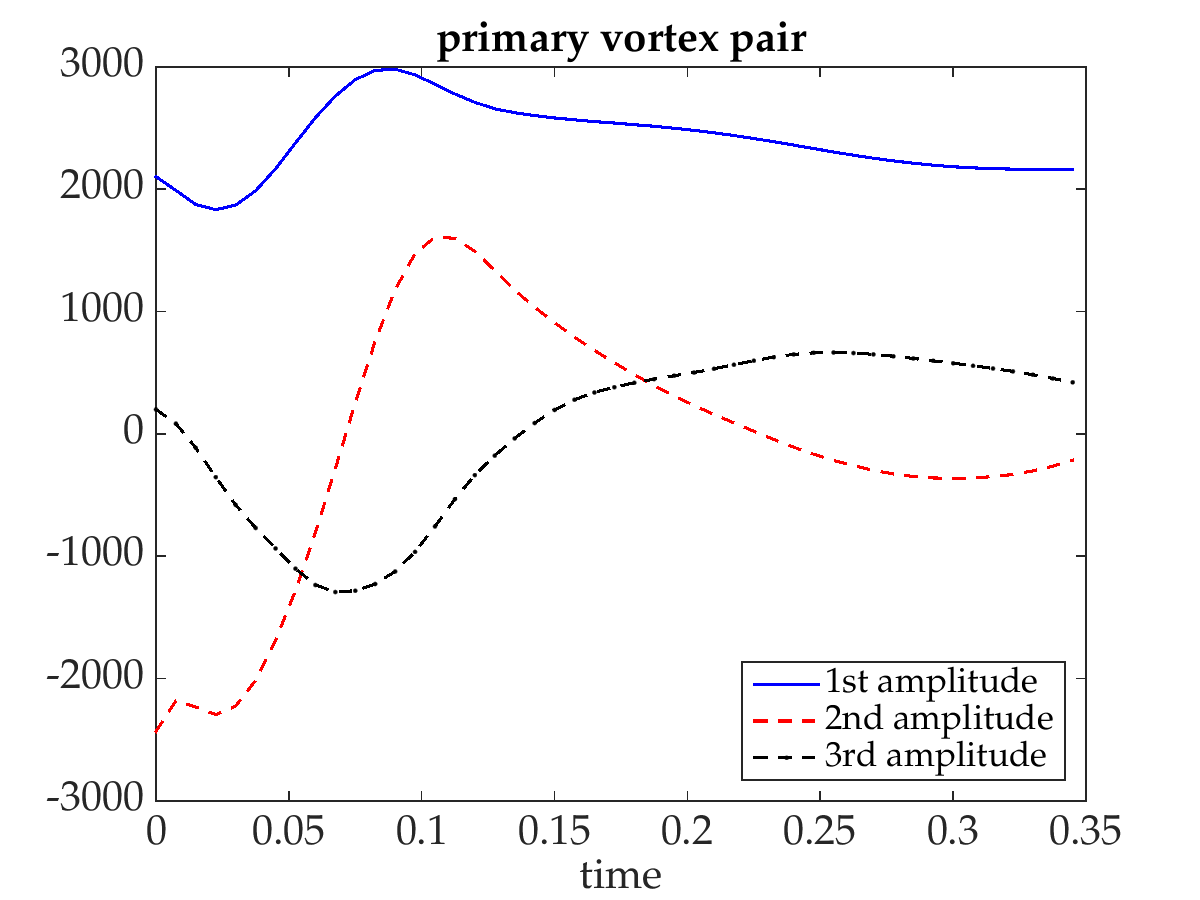}
\includegraphics[width=0.32\linewidth,clip=true, viewport = 20 0 560 420]{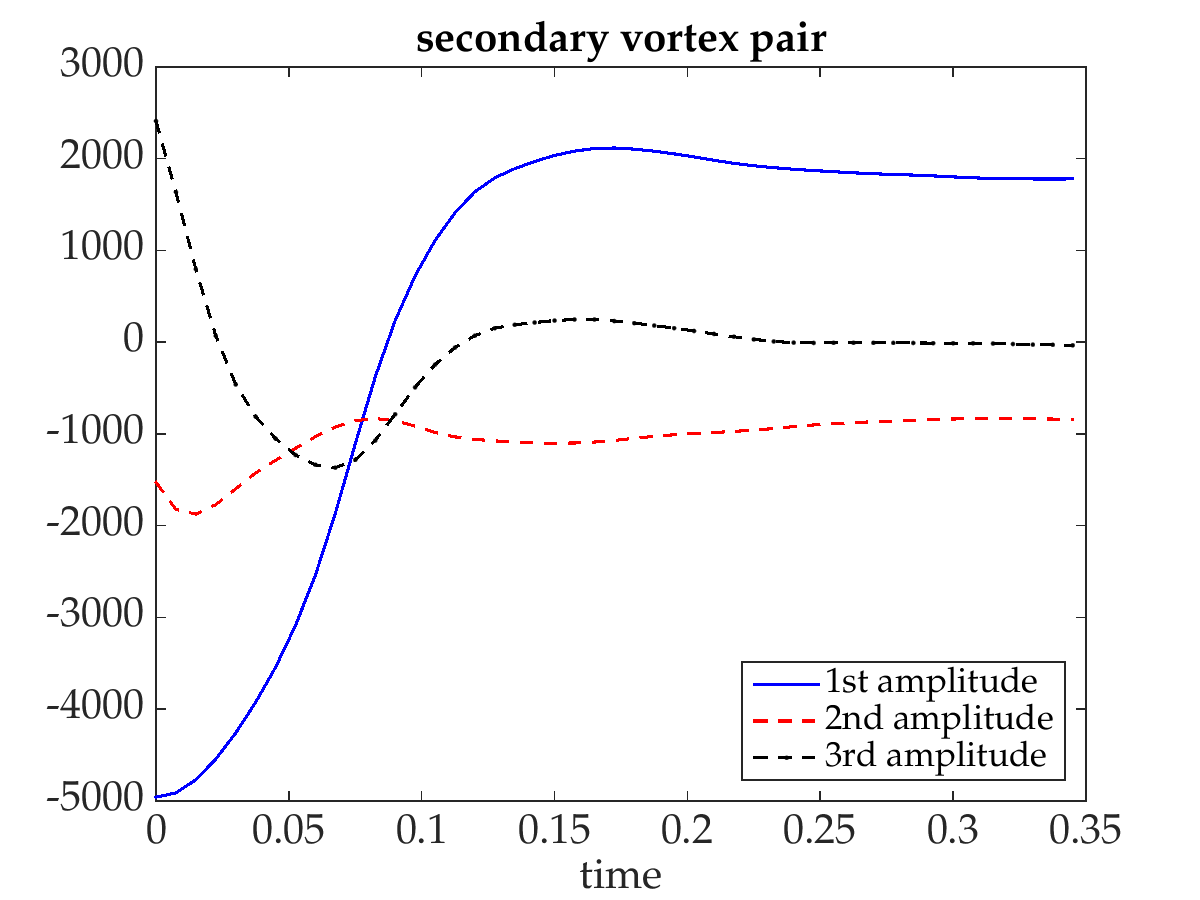}
\includegraphics[width=0.34\linewidth,clip=true, viewport = 10 0 560 420 ]{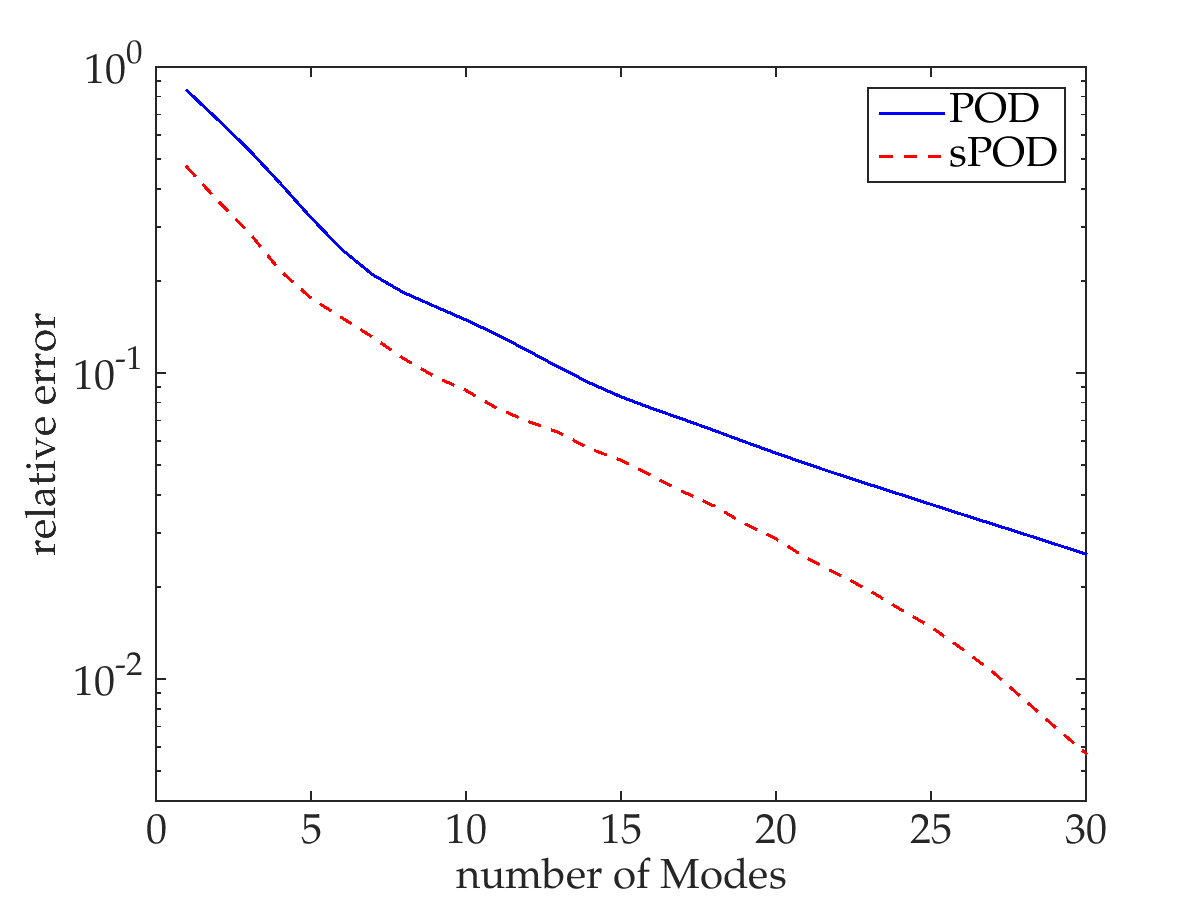}
\caption{2D vortex pair: Left, middle: Time amplitudes for the first three modes for each transported vortex pair.
Right: Mean error over number of modes (POD vs sPOD) }
\label{fig:2D_timeAmplitudes}
\end{figure}

A decent approximation can be obtained by choosing 14  modes for the primary and 10 for the secondary  vortex pair, so in total 24 modes.
The number of modes is selected by the first method described after the sPOD algorithm in section \ref{sec:multipleTransports}.
A decomposition with a large number of modes (here 15 per frame) is calculated and the modes connected with the largest singular values from both frames are selected.  
The decay of the error with this procedure is seen in Fig.\ \ref{fig:2D_timeAmplitudes}, right.
Moreover, the chosen number of modes varies during the sPOD algorithm.
More precisely, we start with one mode per frame and add one mode to each frame in each iteration of the sPOD algorithm until the specified mode number, here 15 per frame, is reached.
The resulting convergence of the sPOD for the two-dimensional vortex pair problem is depicted in Fig.\ \ref{fig:vortexPair_scan}, right. 
In comparison to the pressure pulse example, cf. Fig.\ \ref{fig_sPOD1D}, left, the convergence is worse but still the relative mean error decreases to a value of slightly less than $1\%$.
Due to the used constant shift velocities and the complexity of the transport phenomenon, a much smaller error cannot be expected.
In the insets of Fig.\ \ref{fig:vortexPair_scan}, right, two further convergence criteria are investigated which are introduced in the following.

Two important properties of the SVD are the orthogonality of the approximation to the residual and that the smallest singular value of the approximation is not smaller than the biggest singular value of the residual. These two properties are essential for the optimality of the POD. We formulate their analogues for the case of multiple moving frames:
\begin{equation}
\label{defConvergence}
   \left\langle  \tilde u^k_l  ( \tilde v^k_l)^T , T^{-c_k}( R)\right\rangle_F     = 0\qquad
   \mathrm{and}\qquad
   \min(\tilde \Sigma^k)\geq \max(\Sigma_r^k)
\end{equation}
for all $k=1,\ldots,N_s$. Here $R=X-\tilde X$ is the residual between the snapshot matrix $X$ and the sPOD approximation $\tilde X = \sum_{k=1}^{N_s}T^{c_k} (\tilde U^k \tilde \Sigma^k (\tilde V^k)^T)$ as output from the sPOD algorithm (see section \ref{sec:multipleTransports}). Furthermore, $\tilde u^k_l$ and $\tilde v^k_l$ are the columns of $\tilde U^k$ and $\tilde V^k$, respectively. The matrix $\Sigma_r^k$ is diagonal containing the singular values of $T^{-c_k}( R)$ and $\langle\cdot,\cdot\rangle_F$ denotes the Frobenius inner product.

The first inset in Fig.\ \ref{fig:vortexPair_scan}, right, shows the ratio $\max(\Sigma_r^k)/\min(\tilde \Sigma^k)$ plotted over the iterations for both the primary and the secondary vortex pair.
The criterion $\min(\tilde \Sigma^k)\geq \max(\Sigma_r^k)$ is satisfied from iterations 9 and 12 on, respectively.
In the second inset, the iterations of the maximum relative orthogonality error
\begin{equation*}
\max_l\left( \frac{\left| \left\langle  \tilde u^k_l  \left( \tilde v^k_l\right)^T , T^{-c_k}\left( R\right)\right\rangle_F  \right|}{  \lVert \tilde u^k_l  \left( \tilde v^k_l\right)^T\rVert_F\, \lVert T^{-c_k}\left(R\right) \rVert_F}\right)
\end{equation*}
are shown.
Here, the decay is not as clear as for the other criteria which is due to a dominance by modes associated with a small weight, while the modes associated with larger singular values have a much lower orthogonality error. After all, the decay of the mean error appears to be the most appropriate criterion for describing the convergence.

\begin{figure}
\centering
\includegraphics[width=\linewidth]{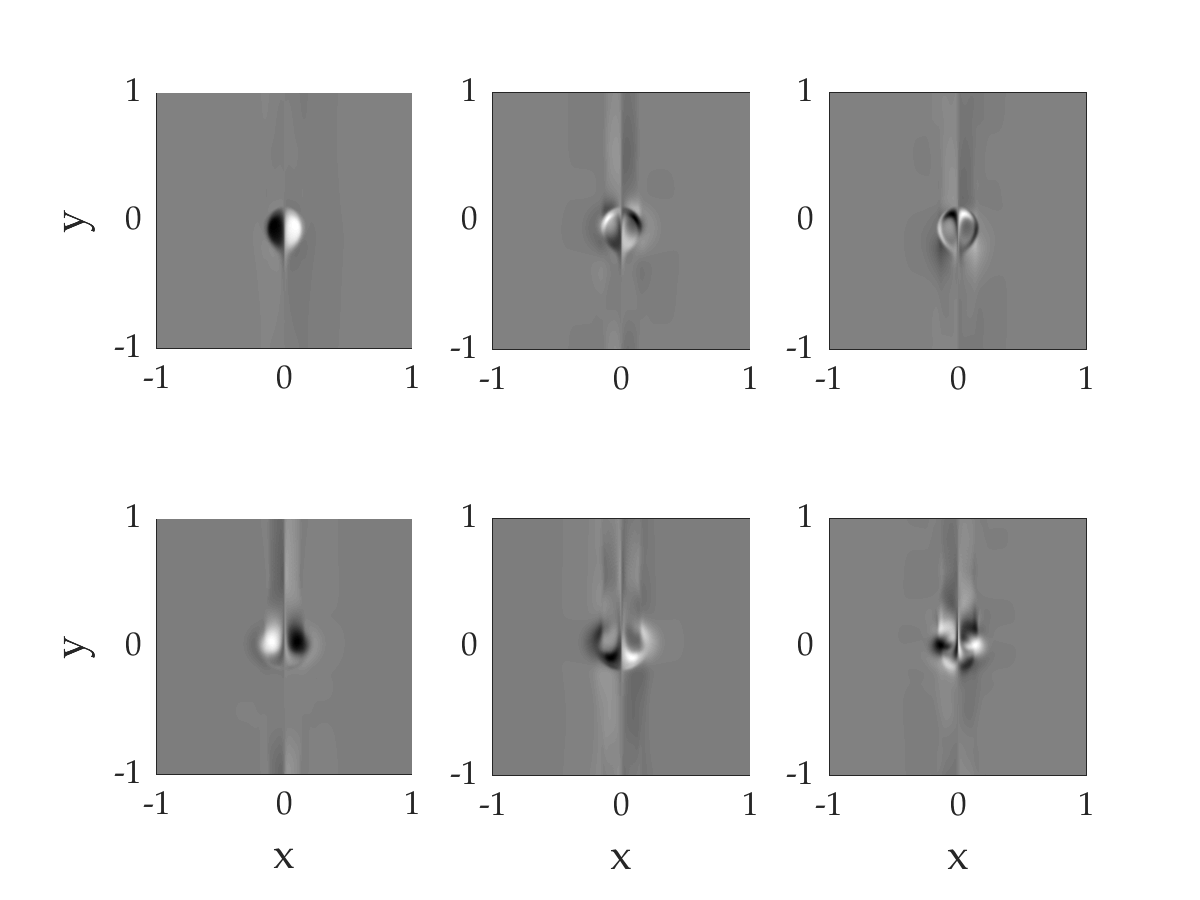}
\caption{2D vortex pair: First three sPOD modes for each transported vortex pair (top: Primary vortex pair, bottom: Secondary vortex pair)}
\label{fig:2D_spatialModes}
\end{figure}

%
The first three modes obtained by the sPOD for each transported vortex pair are shown in Fig.\ \ref{fig:2D_spatialModes}.
The depicted modes correspond to the respective co-moving frames, i.\,e., $\tilde u^k_l$.
For both velocity frames the respective first mode obviously describes the corresponding vortex pair.
The  higher modes seem to fulfill  different roles in both cases.
When considering the dynamical numerical solution, we observe that the secondary, down-going pair strongly changes in time; a rotating vorticity distribution is visible at creation.
It looks very much like the second mode, thus, the higher modes represent the change of shape of the vortex pair.
On the other hand, the numerical simulation results reveal that the shape of the primary, up-going vortex pair changes less in time.
Right after the separation the sharp vortex structure emerges.
The higher modes of the primary vortex pair have a strong weight at the front and rear edge of the pair.
These structures seem to account for a velocity correction, since they induce a translational shift of the primary vortex pair. The velocity correction is caused by the non-constant propagation velocity of the vortex pair.

The corresponding time amplitudes of the spatial modes are depicted in Fig.\ \ref{fig:2D_timeAmplitudes}.
The amplitudes are given by the right singular vectors $\tilde v_l^k $, describing the temporal development of the respective mode $\tilde u^k_l$, multiplied by the corresponding singular value $\tilde \sigma^k_l$ .
For large values of $t$, the first mode dominates the behavior of the primary vortex pair, while for the secondary pair the first two modes have a nearly similar weight. Thus, the first and the second mode of the secondary vortex pair both contribute to its shape. This can be comprehended by considering Fig.\ \ref{fig:2D_spatialModes}, where it can be seen that the vortex pair appearing in the first mode seems to fit into the cavity formed by the structures of the second mode. This agrees with the observation that the secondary vortex pair depicted in Fig.\ \ref{fig:vortexPair_FOM_Solution} is wider than the vortex pair shown in the first mode, cf. Fig.\ \ref{fig:2D_spatialModes}, bottom left.
The separation of the two vortex pairs is completed at about $t>0.15$s. 
The interplay between the modes of the two vortex pairs is complicated up to that point, but afterwards a simpler picture emerges.
For the primary vortex pair, the absolute value of the amplitude of the dominating first mode decays gradually in time, which reflects the reduction of the vortex strength by viscosity.
The influence of the second mode becomes stronger from roughly $t= 0.2$ on, i.\,e., after the full separation.
It leads to a reduction of the vortex pair transport velocity.
Also, the third mode contributes to this velocity decrease.
 In contrast, the initially negative time amplitude of the first mode of the secondary vortex pair increases and changes the sign over time.
 After separation all three amplitudes are nearly constant; this corresponds to the plausible observation that the friction has less impact on the secondary than on the primary vortex pair due to the lower velocity gradient.

While the modes give  a clear and intuitive description of the flow dynamics, still some unexpected structures are visible.
Namely,  one can also see some stripe structures in $y$-direction, especially at the first mode of the secondary vortex pair, cf. Fig.\ \ref{fig:2D_spatialModes}. 
Our first analysis indicates that this is {\it not} due to an incomplete separation, but rather caused by a non-uniqueness of the decomposition. This non-uniqueness can be comprehended by the following consideration.
If we add a constant offset to one of the transported quantities, the combined sPOD approximation remains unchanged if the negative counterpart offset is added to the other transported quantity. Note that offsets are shift-invariant and, hence, well represented in all velocity reference frames.
The occurring stripes seem to be localized variants of these offsets along the $y$-direction.
These structures may be removed by an additional, physically motivated constraint. 
 Nevertheless, the main structure of the vortex pairs is already captured quite well by the respective first modes.

\begin{figure}
\centering
\includegraphics[clip=true, viewport =40 130 530 290, width=\linewidth]{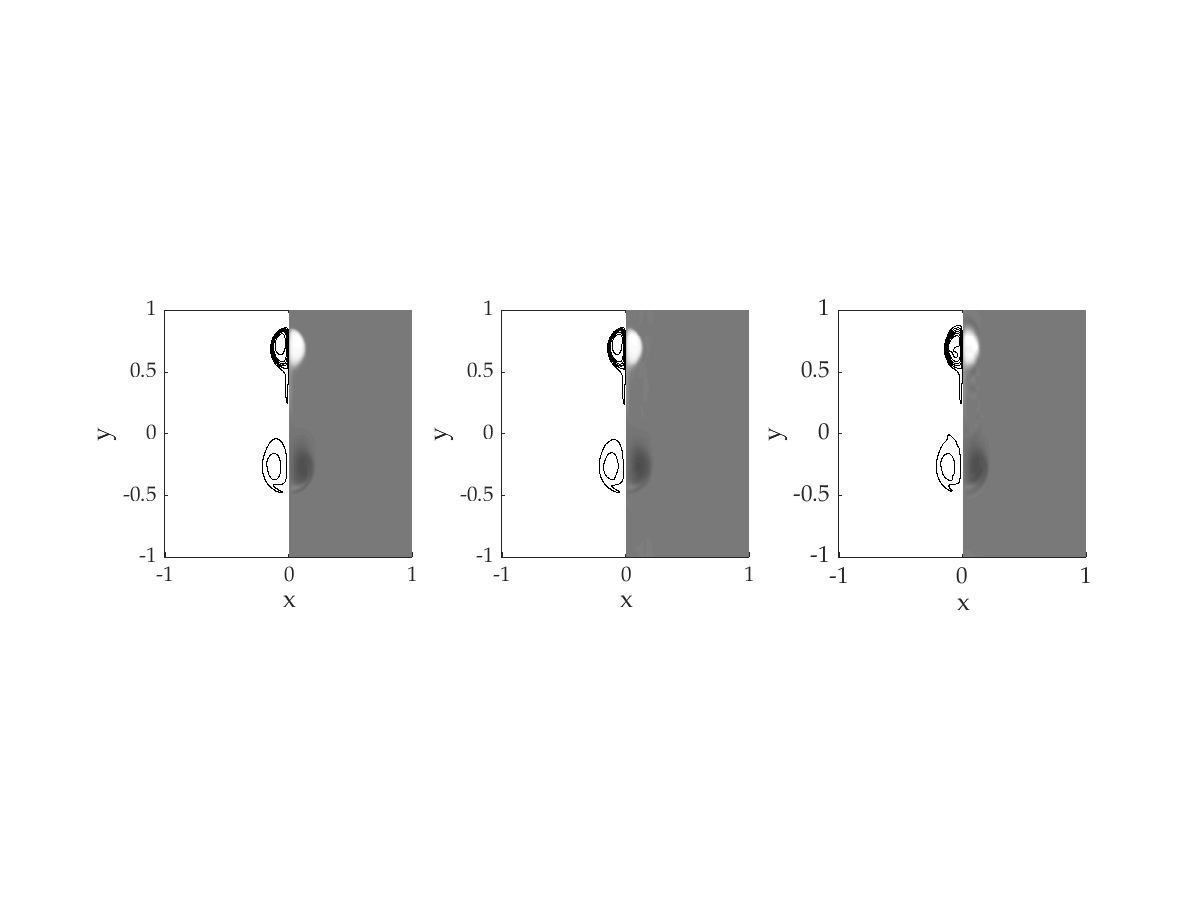}
\caption{ 2D vortex pair: Comparison of full-order solution with sPOD and POD approximation at $t=0.27s$ (left: Full-order solution, middle: sPOD approximation (9 modes for the primary, 6 for the secondary pair), right: POD approximation (15 modes), left half of each plot: Vorticity contour lines at $-270,-250,\ldots,270$, right half of each plot: Pseudocolor plot)}
\label{fig:2D_shPOD_vs_POD}
\end{figure}

The spatial modes and their time amplitudes both provide physically meaningful insights into the dynamics of the vortex movements.
Finally, the sPOD is to be evaluated quantitatively by comparing it to the full-order solution and the POD approximation. To this end, Fig.\ \ref{fig:2D_shPOD_vs_POD} shows the full-order solution, the sPOD approximation ($9 + 6 $ modes), and the POD approximation ($15$ modes), each for a constant time $t=0.27s$.
It is hard to see a difference between the sPOD approximation and the full-order solution, since they match almost exactly. In contrast, even though the POD approximation captures the main structures of the vortex pairs quite well, it reveals some unphysical structures in the inside as well as a blurred surface of the primary vortex pair (top) and, furthermore, some spurious edges in the secondary vortex pair (bottom). The superiority of the sPOD algorithm can also be comprehended by considering Fig.\ \ref{fig:2D_timeAmplitudes}, right,   
 where the mean error is plotted over the number of modes for the POD and the sPOD.
The sPOD mean error is significantly smaller than the POD mean error for all tested numbers of modes although the difference is not as striking as in the wave equation example, cf. section \ref{sec:multipleTransports}. This is not surprising, since the 1D example is ideal for the reduction by the sPOD while the vortex pair problem provides some additional challenges. 
Further, the small transport distance compared with the relatively large interaction time at the beginning reduces the difference between the methods. 

It should be noted that very high numbers of modes may be problematic for the described method. 
The  separation of velocity frames builds on low rank approximations so that the algorithm is expected to become ineffective for a large number of modes.

All in all, the sPOD method performs very well despite the non-trivial and non-constant velocities of the two-dimensional vortex pair problem
and further delivers modes which allow a physically meaningful interpretation.

\section{Summary and outlook\label{sec:conclusion}}

We have presented the shifted proper orthogonal decomposition (sPOD) as a new model reduction approach. It generalizes the common POD by allowing for time-dependent shifts of the modes. More precisely, given a snapshot matrix $X$ and transport velocities $c_1,\ldots,c_{N_s}$, the sPOD constructs an approximate decomposition of the form

\begin{equation}
\label{eq:sPODdecomposition}
X\approx \sum_kT^{c_k}  \left(\tilde U^k \tilde S^k (\tilde V^k)^T\right)
\end{equation}
with low rank matrices $\tilde U^k \tilde S^k (\tilde V^k)^T$ and shift operators $T^{c_k}$ which shift every column of their respective arguments by an amount defined by $c_k$. Due to the shifts, the sPOD is especially well-suited for transport-dominated phenomena. Even problems with multiple shifts can be successfully separated and decomposed based on the iterative sPOD algorithm which is based on truncated SVDs of shifted snapshot matrices and shifted residual matrices. The algorithm can easily be generalized to non-constant velocities as demonstrated in the example of two crossing shocks.

For the case $N_s=1$ and $c_1=0$ in \eqref{eqn:multi_struc}, the common POD is recovered and hence the sPOD delivers good low-dimensional representations for problems where the POD is applied successfully.
Moreover, for the linear wave equation we have demonstrated that the sPOD is able to describe transport phenomena exhibiting sharp moving structures with the minimal number of modes and recovering the analytic solution up to machine precision.
In contrast, for these examples the common POD performs poorly due to the slow singular value decay.
In general, for examples where the snapshot matrix can be perfectly described by a decomposition of the form \eqref{eq:sPODdecomposition}, the analytic solution is a fix point of the sPOD iteration.
However, in general there may be more than one fix point and it is a priori unclear to which one the sPOD iteration converges.
This non-uniqueness has been investigated numerically for the wave equation with two transported quantities by adding a redundant velocity $c_3=0$.
Still the sPOD iteration reduces the error up to any given achievable tolerance although it does not detect that there is a redundancy in the number of modes.

We have tested the sPOD algorithm also for more complex examples, where the analytic solution is unknown.
For the example of two crossing shock waves in 1D and for the vortex pair in 2D, the sPOD shows a better performance than the POD in number of modes needed to ensure a certain accuracy.
Moreover, in the co-moving frames the sPOD modes lend themselves a clear and intuitive description as opposed to the POD modes in the lab frame.

The application of the  sPOD depends on the transport velocities $c_k$.
We have discussed methods for determining good values for the $c_k$ in cases where they are not known from physical considerations, namely for instance by peak or front tracking or by a singular value maximization.
In these cases the sPOD can be computed purely based on snapshot data.
We have also investigated numerically the case when the velocities are slightly under- or overestimated.
For the linear wave equation example with two moving quantities, the sPOD needs a few more modes (depending on how strong the velocities are disturbed) but the singular value decay is still much steeper compared with the POD.
All in all for practical applications, the sPOD qualifies as an alternative to the POD for transport-dominated phenomena yielding a stronger singular value decay and physically more intuitive modes.

The results obtained so far motivate for further development and generalization of the sPOD.
In order to ensure local optimality, the decomposition \eqref{eq:sPODdecomposition} may be directly reformulated in terms of an optimization problem.
This is likely to be more expensive than the algorithm proposed here but in return guarantees local optimality.
Moreover, an extension to non-periodic boundaries is necessary to be able to handle more realistic settings, e.\,g., with reflecting and non-reflecting boundaries.
In addition, the considerations made so far only apply to the mode decomposition of the snapshot matrix.
How the sPOD modes can be used in order to obtain a dynamic reduced-order model with an efficient offline/online decomposition needs to be investigated in the future.
\\

\noindent\textbf{Acknowledgements.} The authors gratefully acknowledge the support by the Deut\-sche Forschungsgemeinschaft (DFG) as part of the collaborative research center SFB 1029 \textit{Substantial efficiency increase in gas turbines through direct use of coupled unsteady combustion and flow dynamics}, project A02 \textit{Development of reduced-order models for pulsed combustion}.

\bibliographystyle{plain}
\bibliography{Refs}
\end{document}